\documentclass[11pt]{article}
\usepackage{latexsym}
\usepackage{graphicx} 
\usepackage[all]{xypic}
\usepackage{color} 
\usepackage{amsmath, amsthm, amssymb}
\usepackage{enumerate}
\usepackage{float}
\usepackage{chemarr}
\usepackage{url}
\usepackage{stackrel}
\usepackage[small,bf]{caption}
\usepackage{wrapfig}

\newtheorem{theorem}{Theorem}[section]
\newtheorem{corollary}[theorem]{Corollary}

\newtheorem{lemma}[theorem]{Lemma}
\newtheorem{example}[theorem]{Example}

\newtheorem{definition}[theorem]{Definition}

\newtheorem{algorithm}[theorem]{Algorithm}
\newtheorem{remark}[theorem]{Remark}

\newtheorem{question}[theorem]{Question}
\newtheorem{assumption}[theorem]{Assumption}

\newcommand{\be}{\begin{equation}}
\newcommand{\ee}{\end{equation}}
\newcommand{\bd}{\begin{displaymath}}
\newcommand{\ed}{\end{displaymath}}
\newcommand{\beal}{\begin{align}}
\newcommand{\enal}{\end{align}}
\newcommand{\been}{\begin{enumerate}}
\newcommand{\enen}{\end{enumerate}}
\newcommand{\beit}{\begin{itemize}}
\newcommand{\enit}{\end{itemize}}

\newcommand{\CRN}{chemical reaction network }
\newcommand{\CRNNoSpace}{chemical reaction network}
\newcommand{\CRNs}{chemical reaction networks }

\newcommand{\CRS}{chemical reaction system }

\newcommand{\famKins}{\mathcal{K}}
\newcommand{\ma}{mass-action }

\newcommand{\op}{\operatorname}

\newcommand{\Net}{\mathfrak{G}}
\newcommand{\Sbar}{S} % no longer {\overline{S}}
\newcommand{\invtPoly}{\mathcal{P}} %stoic. c. class

\newcommand{\xOne}{\widetilde x ^*}
\newcommand{\xTwo}{\widetilde x ^{**}}

\newcommand{\SigSet}{\Sigma}

\newcommand{\SSS}{\mathcal S}
\newcommand{\CC}{\mathcal C}
\newcommand{\RR}{\mathcal R}

\providecommand{\abs}[1]{\lvert#1\rvert}

\newcommand{\fd}{fully open } %using language of Craciun & Feinberg (semiopen)

\newcommand{\ra}{\rightarrow}
\newcommand{\lra}{\leftrightarrows}
\newcommand{\la}{\leftarrow}

\newcommand{\drawBox}{
	\fbox{\begin{minipage}{0.1in}\hspace{0.01in}\vspace{0.05in}\end{minipage}} }

\newcommand{\R}{\mathbb{R}}
\newcommand{\Rplus}{\mathbb{R}_{>0}}
\newcommand{\Rnn}{\mathbb{R}_{\geq 0}} % 'non-negative'

\newcommand{\Z}{\mathbb{Z}}

\newcommand{\ep}{\epsilon}

\newcommand{\im}{\operatorname{Im}}

\begin{document}

\title{Atoms of multistationarity in chemical reaction networks}
%\title{Subnetwork analysis for establishing multistationarity in chemical reaction networks}
\author{Badal Joshi and Anne Shiu} 

\maketitle

% ABSTRACT
\begin{abstract}
Chemical reaction systems are dynamical systems that arise in chemical engineering and systems biology.  In this work, we consider the question of whether the minimal (in a precise sense) multistationary chemical reaction networks, which we propose to call `atoms of multistationarity,' characterize the entire set of multistationary networks.  Our main result states that the answer to this question is `yes' in the context of \fd continuous-flow stirred-tank reactors (CFSTRs), which are networks in which all chemical species take part in the inflow and outflow.  In order to prove this result, we show that if a subnetwork admits multiple steady states, then these steady states can be lifted to a larger network, provided that the two networks share the same stoichiometric subspace.  We also prove an analogous result when a smaller network is obtained from a larger network by `removing species.'  Our results provide the mathematical foundation for a technique used by Siegal-Gaskins {\em et al.}\ of establishing bistability by way of `network ancestry.' 
Additionally, our work provides sufficient conditions for establishing multistationarity by way of atoms and moreover reduces the problem of classifying multistationary CFSTRs to that of cataloging atoms of multistationarity.
%(The main application of our work is for enumerating and classifying small multistationary CFSTRs by way of the atoms.)  
As an application, we enumerate and classify all 386 bimolecular and reversible two-reaction networks.  Of these, exactly 35 admit multiple positive steady states.  Moreover, each admits a unique minimal multistationary subnetwork, and these subnetworks form a poset (with respect to the relation of `removing species') which has 11 minimal elements (the atoms of multistationarity).
\end{abstract}

%-------------------------------------
% FIGURE FOR INTRO - part of the figure of 35 minimal subnetworks:
% Poset with respect to `embedded' relation
\begin{wrapfigure}{r}{0.5\textwidth}
  \vspace{-70pt}
  \begin{center}
    \includegraphics[scale=0.58]{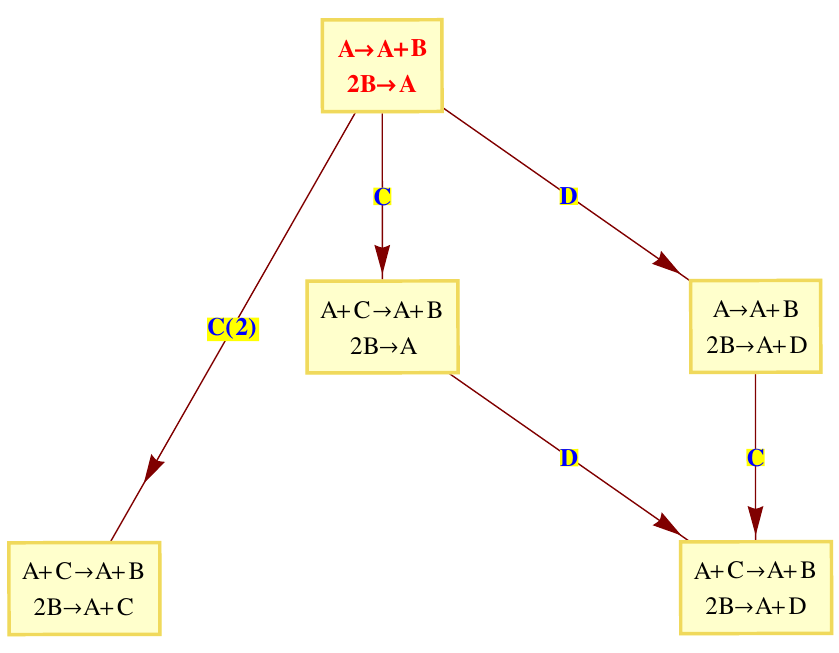}
  \end{center}
    \vspace{-5pt}
  \caption{We propose to call $2B \ra A \ra A+B$ an `atom of multistationarity'; see Section~\ref{sec:atom}. \label{fig:5nets}}
\vspace{-10pt}
\end{wrapfigure}

%  because it `explains' the multistationarity of the four `larger' networks depicted in this figure, which is a portion of Figure~\ref{fig:35}.
%----------------------------------
% SECTION: introduction
%----------------------------------
\section{Introduction}

This work concerns an important class of dynamical systems arising in chemical engineering and systems biology, namely, chemical reaction systems.  As bistable chemical systems are thought to be the underpinnings of biochemical switches, a key question is to determine which systems admit multiple steady states.   In this work, we consider the question of whether the minimal (in a precise sense) networks, which we propose to call `atoms of multistationarity,' characterize the entire set of multistationary networks.  We prove that such atoms do characterize multistationarity for the case of \fd continuous-flow stirred-tank reactors (CFSTRs), which are networks in which all chemical species take part in the inflow and outflow (see Definition~\ref{def:reversible_and_TM}).  For instance, the five networks depicted in Figure~\ref{fig:5nets} are multistationary in the CFSTR setting, but only one is minimal with respect to `removing species' (see Theorem~\ref{thm:enum}).  Following other analyses of small networks~\cite{EnzymeShar,Smallest,SmallGRN,Emergence,ThomasSmallest,smallestHopf}, our main application is to those CFSTRs containing two non-flow reactions.

Chemical reaction systems are nonlinear and parametrized by unknown reaction rate constants.  Thus, determining whether a chemical reaction network admits multiple steady states is difficult: for instance, in the mass-action kinetics setting, it requires determining existence of multiple positive solutions to a system of polynomials with unknown coefficients. However, various criteria have been developed that often can answer this question.  For instance, the Deficiency, Advanced Deficiency, and Higher Deficiency Theories developed by Ellison, Feinberg, Horn, Jackson, and Ji in many cases can affirm that a network admits multiple steady states or can rule out the possibility \cite{EllisonThesis,FeinDefZeroOne,HornJackson72,JiThesis}.  Similarly, the Jacobian Criterion and the more general injectivity test developed by Craciun and Feinberg can preclude multiple steady states~\cite{ME1,ME_entrapped,ME2,ME3,JiThesis,Simplifying}.  These results have been implemented in the CRN Toolbox, freely available computer software developed by Feinberg and improved by Ellison and Ji~\cite{Toolbox}.  Related software programs include BioNetX~\cite{BioNetX,Pantea_comput} and Chemical Reaction Network Software for Mathematica~\cite{CRNsoftware}.  

For systems for which the above software approaches are inconclusive, Conradi~{\em et~al.}\ advocate an approach which first determines whether certain subnetworks admit multiple positive steady states, and if so, tests whether these instances can be lifted to the original network~\cite{Conradi_subnetwork}.  Here, we too examine the topic of lifting multistationarity from a subnetwork to an overall network: Theorem~\ref{thm:subnetwork}, states that this can be accomplished as long as the steady states of interest are nondegenerate and the two networks share the same stoichiometric subspace.  This result and its proof extend Theorem~2 in work of Craciun and Feinberg~\cite{ME_entrapped}.  An important consequence of our theorem is that it provides the mathematical foundation for the technique of Siegal-Gaskins {\em et al.}\ which establishes bistability by way of `network ancestry' (see Remark~\ref{rmk:justify}); their method was applied to a large class of simple gene regulatory networks~\cite{Emergence}.
%A criterion for multistationarity for systems whose steady state locus is defined by binomials is given by P\'erez Mill\'an~{\em et~al.}\ in~\cite{TSS}.  Finally, a degree theory method which can establish multiple steady states can be found in work of Craciun, Helton, and Williams~\cite{CHW08}.

A \fd {\em continuous-flow stirred-tank reactor } (CFSTR) is a network in which all chemical species enter the system at constant rates and are removed at rates proportional to their concentrations (see Definition~\ref{def:reversible_and_TM}).  
In the setting of these systems, we extend our results beyond subnetworks to `embedded networks' which are obtained by removing species as well as reactions from a network (Definition~\ref{def:reduced_network}).  Corollary~\ref{cor:liftEmbCFSTR} states that if an embedded CFSTR of a CFSTR is multistationary then so is the CFSTR itself.  Therefore, the set of multistationary CFSTRs is characterized by its minimal elements (with respect to the embedded network relation), and we pose the challenge of characterizing these atoms.  In this work, we focus on cataloging the smallest atoms.  % A consequence of Theorem~\ref{thm:subnetwork} is that a CFSTR admits multiple nondegenerate positive steady states if some sub-CFSTR does, with no requirement that they share the same stoichiometric subspace (Corollary~\ref{cor:liftCFSTR}).

Recent work of the first author presented a simple characterization of the one-reaction \fd CFSTRs that admit multiple steady states in the mass-action kinetics setting (Theorem~\ref{thm:1rxn})~\cite{1rxn}.  Here we consider the bimolecular two-reaction CFSTRs; a network is `bimolecular' if all of its chemical complexes contain at most two molecules.  We enumerate all 386 reversible such networks.  Of these, exactly 35 admit multiple positive steady states.  Moreover, each admits a unique minimal multistationary sub-CFSTR, and these subnetworks form a poset with respect to the embedded network relation that has 11 minimal elements (Theorem~\ref{thm:enum}).  These 11 networks are precisely the CFSTR atoms of multistationarity in the bimolecular two-reaction setting.
% shown to exist in Theorem~\ref{thm:emb}. %Conjecture~\ref{conj:atom} proposes that the set of such minimal networks characterize the multistationary CFSTRs, beyond the bimolecular two-reaction setting.  
Note that a similar enumeration of small bimolecular networks was undertaken by Deckard, Bergmann, and Sauro~\cite{Deckard}, from which Pantea and Craciun sampled networks to compute the fraction of such networks that pass the Jacobian Criterion \cite{Pantea_comput}.
%Recall that one test for precluding multistationarity, which is implemented in some of the above-mentioned software programs, is the Jacobian Criterion.  The criterion, which is due to Craciun and Feinberg, gives sufficient conditions for precluding multiple steady states~\cite{BanajiCraciun2009, BanajiCraciun2010, CPS, ME1,ME_entrapped,ME2,ME3}.  A chemical reaction network is said to pass the Jacobian Criterion (or to be `injective') if all terms in the determinant expansion of its parametrized Jacobian matrix have the same sign. One aim of our work is to demonstrate the the Jacobian Criterion rules out many small networks from admitting multiple steady states.  By applying the methods for simplifying the Jacobian Criterion introduced by the authors in~\cite{Simplifying}, we find that only 55 fail the Jacobian Criterion, of which, as mentioned above, 35 admit multiple steady states.

% REMARK: CFSTR
% Note that although our results are stated for continuous-flow stirred-tank reactors (CFSTRs), when multiple steady states are precluded by the Jacobian Criterion, then we can still rule out multiple steady states when some inflows or outflows are removed, if the underlying non-flow network of interest satisfies some mild hypotheses (such as weakly-reversibility; see Definition~\ref{def:reversible_and_TM}) \cite[Corollary~8.3]{ME3}.

%:OUTLINE of article
This article is organized as follows.  
% Section 2
Section~\ref{sec:intro_CRN} introduces chemical reaction systems.  
% Section 3
Our main result for lifting multiple steady states from subnetworks, Theorem~\ref{thm:subnetwork}, appears in Section~\ref{sec:lift}.
% Section 4
In Section~\ref{sec:lift_emb}, this result is extended in the case of \fd CFSTRs: Theorem~\ref{thm:emb} implies that steady states from embedded CFSTRs can be lifted as well. 
% Section 5
Section~\ref{sec:atom} introduces `atoms of multistationarity.'  
% Sections 6-7
Section~\ref{sec:enum} describes our approach to enumerating bimolecular two-reaction networks (Algorithm~\ref{algor:enum}), and Section~\ref{sec:class} determines which such networks are multistationary in the mass-action kinetics setting (Theorem~\ref{thm:enum}) and displays the resulting atoms of multistationarity (Figure~\ref{fig:35}).

%----------------------------------
% SECTION: Chemical reaction network theory
%----------------------------------
\section{Chemical reaction network theory} \label{sec:intro_CRN}
In this section we review the standard notation and recall the classification of one-reaction CFSTRs.  % in Section~\ref{sec:1Rxn}.
% and the use of the Jacobian Criterion to rule out multistationarity in Section~\ref{sec:JC}.  

% SUBSECTION: CRNs
\subsection{Chemical reaction networks}
We begin with an example of a {\em chemical reaction}:
\begin{align*}
  2X_{1}+X_{2} ~\rightarrow~ X_{3}~.
\end{align*}
Each $X_{i}$ is called a chemical {\em species}, and $2X_{1}+X_{2}$ and
$X_{3}$ are called chemical {\em complexes.}  Assigning the {\em
  reactant} complex $2X_{1}+X_{2}$ to the vector $y =
(2,1,0)$ and the {\em product} complex $X_{3}$ to the vector
$y'=(0,0,1)$, we can write the reaction as $ y \rightarrow y'$. In
general we let $s$ denote the total number of species $X_{i}$, and we
consider a set of $r$ reactions, each denoted by
\begin{align*}
  y_{k} \rightarrow y_{k}'~, 
\end{align*} 
for $k \in \{1,2,\dots,r\}$, and $y_k, y_k' \in \Z^s_{\ge 0}$, with $y_k \ne
y_k'$.  We index the entries of a 
complex vector $y_k$ by writing $y_k = \left( y_{k1}, y_{k2}, \dots, y_{ks}\right) \in \Z^s_{\geq 0}$,
and we will call $y_{ki}$ the {\em stoichiometric coefficient} of species $i$ in 
complex $y_k$.  
For ease of notation, when there is no need for
enumeration we typically will drop the subscript $k$ from the notation
for the complexes and reactions.  

% Def: self-catalyzing
%A {\em self-catalyzing reaction} will refer to a reaction, such as $A+B \ra 2A+C$, in which there is a species (here, $A$) that is a {\em self-catalyst}: it appears with a higher stoichiometric coefficient in the product complex than in the reactant complex. 
% DEFINE: CRN
\begin{definition}   \label{def:crn}
  Let $\SSS = \{X_i\}$, $\CC = \{y\},$ and $\RR = \{y \to y'\}$ denote finite 
  sets of species, complexes, and reactions, respectively.  The triple
  $\{\SSS, \CC, \RR \}$ is called a {\em chemical reaction network} if it satisfies the following:
  \been
  	\item for each complex $y \in \CC$, there exists a reaction in $\RR$ for which $y$ is the reactant complex or $y$ is the product complex, and
	\item for each species $X_i \in \SSS$, there exists a complex $y \in \CC$ that contains $X_i$.
  \enen
\end{definition}
%Later we will view $\RR$ as a multiset. 

%We say that a species $X_i$ is a {\em reactant species} (respectively, {\em product species}) if it appears in a reactant (respectively, product) complex of some reaction. 
% Define: decouple
A network {\em decouples} if there exist nonempty subsets $\RR' \subset \RR$ and $\RR'' \subset \RR$ such that $\RR = \RR' ~ \dot \cup ~ \RR''$ and such that the species involved in reactions in $\RR'$ are distinct from those of $\RR''$. %Otherwise, we say that the network is {\em coupled}. 
% Define: subnetwork
We next define a subnetwork and the more general concept of an `embedded' network,'
which was introduced by the authors in \cite[\S 4.2]{Simplifying}.  Informally, a network $N$ is an embedded network of a network $G$ if $N$ may be obtained from $G$ by removing reactions and `removing species.'
% DEFINITION: reduced network
\begin{definition} \label{def:reduced_network}
Let $\Net=\{\SSS,\CC,\RR\}$ be a \CRNNoSpace.
\begin{enumerate}
	\item Consider a subset of the species $S\subset \SSS$, a subset of the complexes $C \subset \CC$, and a subset of the reactions $R\subset \RR$. 
	\beit
		\item 
		The {\em restriction of $R$ to $S$}, denoted by $R|_S$, is the set of reactions obtained by taking the reactions in $R$ and removing all species not in $S$ from the reactant and product complexes. If a trivial reaction (one in which the reactant and product complexes are the same) is obtained in this process, then it is removed.  Also removed are extra copies of repeated reactions.
		\item The {\em restriction of $C$ to $R$}, denoted by $C|_R$, is the set of (reactant and product) complexes of the reactions in $R$.  
		\item The {\em restriction of $S$ to $C$}, denoted by $S|_C$, is the set of species that are in the complexes in $C$. 
	\enit
	
	\item The network obtained from $\Net$ by {\em removing a subset of species} $\{X_i\} \subset \SSS$ is the network 
		$$\left\{\SSS\setminus \{X_i\} ,~\CC|_{\RR|_{\SSS\setminus \{X_i\}}},~
			\RR|_{\SSS\setminus \{X_i\}}  \right\}~.$$
	\item A subset of the reactions $\RR' \subset \RR$ defines the {\em subnetwork} $\{\SSS|_{\CC|_{\RR'}},\CC|_{\RR'},\RR' \}$. %$ where $\CC|_{\RR'}$ denotes the set of complexes that appear in the reactions $\RR'$, and $\SSS|_{\CC_{\RR'}}$ denotes the set of species that appear in those complexes. 
	%The network obtained from $\Net$ by {\em removing a set of reactions} $\{y \ra y'\} \subset \RR$ is the subnetwork 
%		$$\left\{\SSS|_{\CC|_{\RR \setminus \{y \ra y' \}}},~\CC|_{\RR \setminus \{y \ra y' \}},~\RR \setminus \{y \ra y' \} \right\}~.$$
	% DEFINITION: embedded network
	\item
	Let $\Net=\{\SSS,\CC,\RR\}$ be a \CRNNoSpace.
	An {\em embedded network} of $\Net$, which is defined by a subset of the species, $S=\left\{X_{i_1}, X_{i_2},\dots, X_{i_{k}} \right\}\subset \SSS$,
	and 
	a subset of the reactions, $R=\left\{ R_{j_1}, R_{j_2},\dots, R_{j_l} \right\}\subset \RR$, that involve all species of $S$, 
	is the network $(S,\CC|_{R|_S},R|_S)$ consisting of the reactions $R|_S$. % in which all species not in $S$ are removed.
\end{enumerate}
\end{definition}
\noindent
% Similar definition appears in Anderson's latest paper (what he calls $S$-reduced networks), Casian too.
\begin{remark}
We note that a network is also a subnetwork and an embedded network of itself.  In fact, any subnetwork $\{\SSS|_{\CC|_{\RR'}},\CC|_{\RR'},\RR' \}$ is an embedded network, namely the one defined by the subset of species $\SSS|_{\CC|_{\RR'}}$ and the subset of reactions $ \RR'$.   

We also note for readers who are familiar with species-reaction (SR) graphs that the definitions of `subnetwork' and `embedded network' can be interpreted as follows.  Recall that the SR graph of a network consists of species vertices and reaction vertices, with edges arising from reactions in the network; for details, see~\cite{ME2}.  A subnetwork corresponds to the subgraph of the SR graph induced by the full set of species and the subset of reaction vertices arising from reactions in the subnetwork.  As for an embedded network, this arises as the subgraph induced by the corresponding subsets of species and reaction nodes.
\end{remark}
One focus of our work is on CFSTRs, which we now define.
% DEFINITION: weakly-reversible (removed)
% DEFINITION: CFSTR
\begin{definition} \label{def:reversible_and_TM}
\begin{enumerate}
%	\item  A \CRN is {\em weakly-reversible} if each complex of the network is both a reactant complex and a product complex.
%  A reaction $y \to y' \in
%  \RR$ is {\em reversible} if its reverse reaction $y' \to y$ is also a reaction in $\RR$. 
   %A network is said to be {\em reversible} if all of its reactions are reversible. 
	\item 
	A {\em flow reaction} contains only one molecule; such a reaction is either an {\em inflow reaction} $0 \ra X_i$ or an {\em outflow reaction} $X_i \ra 0$.
	\item   % DEFINITION OF CFSTR
  A \CRN is a {\em continuous-flow stirred-tank reactor} {\em (CFSTR)} if it contains all outflow reactions $X_i \ra 0$ (for all $X_i \in \SSS$). A CFSTR is {\em fully open} if it additionally contains all inflow reactions $0 \ra X_i$.  %Note that each reaction 
%network $\Net=\{\SSS, \CC, \RR \}$ is contained in a unique minimal \fd CFSTR which we denote by $\widetilde \Net := \{\SSS,~\CC \cup \SSS \cup \{0\},~\RR \cup \{0 \lra X_i\}_{X_i \in \SSS} \}$.  
A {\em sub-CFSTR} is a subnetwork that is also a CFSTR.
  \end{enumerate}
\end{definition}
We note that Craciun and Feinberg use the term `feed reactions' for inflow reactions and `true reactions' for non-flow reactions.
In chemical engineering, a CFSTR refers to a well-mixed tank in which reactions occur.  An inflow reaction represents the flow of species (at a constant rate) into the tank in which the non-flow reactions take place, and an outflow reaction represents the removal or degradation of a species (at rate proportional to its concentration).

\begin{example} \label{ex:running}
Consider the following \fd CFSTR:
	\begin{align} \label{eq:running}
	% INFLOW REACTIONS
	0 \stackrel[\kappa_2]{\kappa_1}{\rightleftarrows} A \quad \quad 
	0 \stackrel[\kappa_4]{\kappa_3}{\rightleftarrows} B \quad \quad 
	0 \stackrel[\kappa_6]{\kappa_5}{\rightleftarrows} C \quad \quad 
	% NON-flow REACTIONS
	2A \stackrel[\kappa_8]{\kappa_7}{\rightleftarrows} A+B \stackrel[\kappa_{10}]{\kappa_{9} }{\rightleftarrows} A+C~.
	\end{align}
%The total molecularity of species $A$, $B$, and $C$ are 5, 2, and 1, respectively.  
% SUBNETWORK
The following sub-CFSTR arises by removing two reactions: % (and it is the unique minimal multistationary subnetwork): 
	\begin{align} \label{eq:running_sub}
	% INFLOW REACTIONS
	0 \stackrel[\kappa_2=1]{\kappa_1}{\rightleftarrows} A \quad \quad 
	0 \stackrel[\kappa_4=1]{\kappa_3}{\rightleftarrows} B \quad \quad 
	0 \stackrel[\kappa_6=1]{\kappa_5}{\rightleftarrows} C \quad \quad 
	% NON-flow REACTIONS
	2A \stackrel[\kappa_8]{~}{\leftarrow} A+B \stackrel[\kappa_{10}]{~}{\leftarrow} A+C~.
	\end{align}
%EMBEDDED NETWORK
Next, we obtain the following embedded network by removing species $C$: 
	\begin{align} \label{eq:running_emb}
	% INFLOW REACTIONS
	0 \stackrel[\kappa_2=1]{\kappa_1}{\rightleftarrows} A \quad \quad 
	0 \stackrel[\kappa_4=1]{\kappa_3}{\rightleftarrows} B \quad \quad 
%	0 \stackrel[\kappa_6=1]{\kappa_5}{\rightleftarrows} C \quad \quad 
	% NON-flow REACTIONS
	2A \stackrel[\kappa_8]{~}{\leftarrow} A+B \stackrel[\kappa_{10}]{~}{\leftarrow} A~.
	\end{align}
\end{example}

%--------------------------------------
%SUBSECTION: dynamics and steady states
%--------------------------------------
\subsection{Dynamics and steady states}
% Differential equations
The concentration vector 
\begin{align*}
	x(t) = \left(x_1(t), x_2(t), \ldots, x_{|\SSS|}(t) \right)
\end{align*}
 will track the concentration $x_i(t)$ of the $i$-th species at time $t$.
A chemical reaction network defines a dynamical system by way of
a {\em rate function} for each reaction.  In other words, to each reaction $y_k \to y_k'$ we assign a %continuously differentiable 
smooth function $\displaystyle R_k(\cdot) =
R_{y_k \to y_k'}(\cdot)$ that satisfies the following assumption.
% ASSUMPTION: rate functions
\begin{assumption} 
 For $k \in \RR$, $\displaystyle R_k(\cdot) = R_{y_k
    \to y_k'}(\cdot): \R^{|\SSS|}_{\ge 0} \to \R$ satisfies:
  \begin{enumerate}
  \item $R_{y_k \to y_k'}(\cdot)$ depends explicitly upon $x_i$ only
    if $\displaystyle y_{ki} \ne 0$.
  \item $\displaystyle \frac{\partial}{\partial x_i}R_{y_k \to
      y_k'}(x) \ge 0$ for those $x_i$ for which $\displaystyle y_{ki}
    \ne 0$, and equality can hold only if at least one coordinate of $x$ is zero.
%    $x \in \partial \R^{|\SSS]}_{\ge  0}$.
  \item $R_{y_k \to y_k'}(x) = 0$ if $x_i = 0$ for some $i$ with
    $\displaystyle y_{ki} \ne 0$.
  \item If $\displaystyle 1 \le y_{ki} < y_{\ell i}$, then
    $\displaystyle \lim_{x_i \to 0} \frac{R_{\ell}(x)}{R_{k}(x)} = 0$,
    where all other $x_j > 0$ are held fixed in the limit.
  \end{enumerate}
  \label{assump:rates}
\end{assumption}

The final assumption simply states that if the $l$-th reaction demands
strictly more molecules of species $X_i$ as inputs than does the $k$-th
reaction, then the rate of the $l$-th reaction decreases to zero faster
than the $k$-th reaction, as $x_i \to 0$. The functions $R_{k}$ are called the {\em kinetics} of the system.

%-------------
% DEFINITION: CRS
%-------------
\begin{definition} \label{def:CRS}
Consider a \CRN $\left\{ \SSS, \CC, \RR=\{y_k \ra y_k' \} \right\}$ and a choice of kinetics $\{R_k\}$ that satisfy Assumption~\ref{assump:rates}.
	\begin{enumerate}
	% DEFINITION: CRS
	\item The following system of ODEs defines a dynamical system is called a {\em chemical reaction system}:
\begin{equation}
  \dot x(t) \quad = \quad  \sum_{k=1}^{\RR} R_{k}(x(t))(y_k' - y_k) \quad =: \quad f(x(t))~,
  \label{eq:main_general}
\end{equation}
where the second equality is a definition.
	% DEFINITION: stoic. c. class
	\item
	The {\em stoichiometric subspace} of the network is the span of all {\em reaction vectors} $y_k' - y_k$.  We
	will denote this space by $\Sbar$ and its dimension by $\sigma$:
		\begin{equation*}
		\Sbar \quad:=\quad {\rm span} \left\{ y_1' - y_1, ~ y_2'-y_2, ~ \dots ~ , ~ y_{|\RR|}' - y_{|\RR|} \right\}~\subset ~ \R^{|\SSS|} ~.
		\end{equation*}
	Note that~\eqref{eq:main_general} implies that a trajectory $x(t)$ that begins at a positive vector $x(0)=c^0 \in \R^s_{>0}$ 
	remains in the  
	{\em stoichiometric compatibility class}, which we denote by
		\begin{align}\label{eqn:invtPoly}
		\invtPoly \quad := \quad (c^0+\Sbar) \cap \mathbb{R}^{|\SSS|}_{\geq 0}~, 
		\end{align}
	for all positive time; in other words, this set $\invtPoly$ is forward-invariant with respect to~\eqref{eq:main_general}.  Two points in the same stoichiometric compatibility class $\invtPoly$ are said to be {\em stoichiometrically compatible}.
	% DEFINITION: STEADY STATE (and `nondegenerate' and `stable')
	\item
	A concentration vector $\overline x \in \R^{|\SSS|}_{> 0}$ is a (positive) {\em steady state} of the system~\eqref{eq:main_general} if $f(\overline x) = \overline{0}$.  
	A steady state $\overline x$ is {\em nondegenerate} if $\im df (\overline x) = \Sbar $.
	%$\left( \ker {d}f (\overline x) \right) \cap \overline S = \{\overline 0 \}$.  
	(Here, ``$df(\overline x)$'' is the Jacobian matrix of $f$ at $\overline x$: the ${|\SSS|} \times {|\SSS|}$-matrix whose $(i,j)$-th entry is equal to the partial derivative $\frac{\partial f_i}{\partial x_j} (\overline x)$).
	A nondegenerate steady state $\overline x$ is {\em exponentially stable} if each of the $\sigma:=\dim \Sbar$ nonzero eigenvalues of $df(\overline x)$ (viewed over the complex numbers) has negative real part.
	\end{enumerate}
\end{definition}
\noindent
% We use the notation $\Sbar$ for the stoichiometric subspace to match that of Craciun and Feinberg~\cite{ME_entrapped}.
%, who distinguish between the span of all non-flow reaction vectors (which they denote by $S$) and that of all flow and non-flow reactions (which we too denote by $\Sbar$)~\cite{ME_entrapped}.
% In CFSTR case: stoich. c. class is all of R^s
In the case of a CFSTR, the reaction vector for the $i$-th inflow reaction is the $i$-th canonical basis vector of $\R^{|\SSS|}$, so the stoichiometric subspace is $\Sbar = \R^{|\SSS|}$.  It follows that for a CFSTR, the unique stoichiometric compatibility class is the nonnegative orthant: $\invtPoly=\R^{|\SSS|}_{\geq 0}$.

% Mass-action kinetics
An important example of kinetics is {\em mass-action kinetics}; a
chemical reaction system is said to have mass-action kinetics if all rate functions $R_{k}$ 
take the following multiplicative form:
\begin{equation}
  R_{k}(x) \quad = \quad  \kappa_k x_1^{y_{k1}} x_2^{y_{k2}} \cdots x_{|\SSS|}^{y_{k{|\SSS|}}}
 \quad  =: \quad \kappa_k x^{y_k}~,
  \label{eq:massaction}
\end{equation}
for some vector of positive {\em reaction rate constants} 
		$(\kappa_1, \kappa_2, \dots, \kappa_{|\RR|}) \in \Rplus^{|\RR|} $, with the convention that $0^0 = 1$.  % and the final equality is a definition.   
It is easily verified that each $R_k$ defined via
\eqref{eq:massaction} satisfies Assumption \ref{assump:rates}.
Combining \eqref{eq:main_general} and \eqref{eq:massaction} gives the
following system of mass-action ODEs:
\begin{equation}
		  \dot x(t) \quad = \quad \sum_{k=1}^{|\RR | } \kappa_k x(t)^{y_k}(y_k' - y_k) \quad =: \quad f(x(t))~.
  \label{eq:main}
\end{equation}
%Here, we make use of multi-index notation: $x(t)^{y_k}:= x_1(t)^{y_{k1}} x_2(t)^{y_{k2}} \cdots x_s(t)^{y_{ks}} $.  
	%(By convention, $0^0 := 1$).  

In the following example and all others in this work, we will label species by distinct letters such as $A,B,\dots$ rather than $X_1,X_2,\dots$.
% EXAMPLE: (continue the running example)
% mass-action ODEs
\begin{example} \label{ex:diffEq}
We now return to the CFSTR~\eqref{eq:running} in Example~\ref{ex:running}.  
The mass-action differential equations~\eqref{eq:main} for this network are the following:
	\begin{align} \label{eq:ex_ODEs}
	\frac{d{ x_A } }{dt} \quad &= \quad  { \kappa_1 } - \kappa_2 { x_{A} } 
		- \kappa_7 x_A^2 + { \kappa_{8} }{ x_{A}x_{B} } 	\notag
	 \\
	\frac{dx_{B}}{dt} \quad &=\quad	{ \kappa_3 } - \kappa_4 { x_{B} } +
		\kappa_7 x_A^2
		-{ \kappa_{8}}{ x_{A}x_{B}}  
  -  \kappa_9 x_A x_B	+ { \kappa_{10} }{ x_{A}x_{C} }
	\\
	\frac{dx_{C}}{dt} \quad &=\quad  { \kappa_5 } - \kappa_6 { x_{C} }
		+ \kappa_9 x_A x_B
		-{ \kappa_{10} }{ x_{A}x_{C} }				\notag
 ~.
	\end{align}
\end{example}

% Next: family of kinetics
Note that a \CRN gives rise to a family of mass-action kinetics systems parametrized by a choice of one reaction rate constant $\kappa_k \in \Rplus$ for each reaction, and all reactions not in the network can be viewed as having reaction rate constant equal to zero.  We now generalize this concept of a parametrized family for other kinetics.

% DEFINITION: MSS
\begin{definition} \label{def:MSS}
\begin{enumerate}
\item A {\em parametrized family of kinetics} $\famKins$ for \CRNs on $|\SSS|$ species is an assignment to each possible reaction $y_k \to y_k'$ (that involves only species from $\SSS$) a smooth function
	\begin{align*}
	\Rnn \times \Rnn^{\SSS} & \to \R^{\SSS} \\
	(\kappa_k,x) &  \mapsto R_{k}^{\kappa_k} (x)
	\end{align*}
such that
	\begin{itemize}
	\item for $\kappa_k>0$, the function $R_{k}^{\kappa_k} (x)$ is a rate function for the reaction $y_k \to y_k'$ that satisfies Assumption~\ref{assump:rates}, and
	\item when $\kappa_k=0$, then $R_{k}^{\kappa_k} (x)$ is the zero function.
	\end{itemize}
% def: closed under scaling - now commented out
%\item A parametrized family of kinetics $\famKins$ is {\em closed under scaling} if for all reactions $k$ and for all $\kappa_k>0$ and $\mu>0$, the scaled rate function $\mu R_{k}^{\kappa_k} (x)$ is also in the family, i.e., there exists $\kappa_k'>0$ such that $ R_{k}^{\kappa_k} (x) = \mu R_{k}^{\kappa_k} (x)$.
% def: MSS	
\item Let $\Net$ be a \CRNNoSpace, and let $\famKins$ be a parametrized family of kinetics on $|\SSS|$ species.  Then $\Net$ is said to admit {\em multiple $\famKins$ steady states} or is {\em $\famKins$-multistationary} if there exist kinetics $\{R^{\kappa_k}_k\}$ arising from $\famKins$ and a stoichiometric compatibility class $\invtPoly$ such that the resulting system~\eqref{eq:main_general} has two or more positive steady states in $\invtPoly$. 
Moreover, such a network is said to admit {\em bistability} if such steady states can be found that are stable.
\end{enumerate}
\end{definition}
As noted above, an important family of kinetics $\famKins$ is that of mass-action kinetics; in this case, $\Net$ admits {\em multiple mass-action steady states} if there exist rate constants $\kappa_{k} \in \mathbb{R}_{> 0}$ 
and a stoichiometric compatibility class $\invtPoly$ such that the mass-action system~\eqref{eq:main} admits at least two positive steady states in $\invtPoly$.

% EXAMPLE: (return to running example - it admits MSS)
\begin{example} \label{ex:2}
We again consider the CFSTR~\eqref{eq:running} examined in Examples~\ref{ex:running} and~\ref{ex:diffEq}.  Recall that for a CFSTR, the unique stoichiometric compatibility class is the nonnegative orthant: here, $\invtPoly=\R^3_{\geq 0}$.  Therefore, our CFSTR~\eqref{eq:running} admits multiple positive mass-action steady states if and only if there exist reaction rate constants $\kappa_1, \kappa_2, \dots, \kappa_{10} \in \R_{>0}$ such that the differential equations~\eqref{eq:ex_ODEs} have at least two positive steady states.  Indeed, the CRN Toolbox~\cite{Toolbox} determines that when the mass-action system takes the following rate constants:
	\begin{small}
	\begin{align*}
	& \left( \kappa_1,\kappa_2, \kappa_3, \kappa_4, \kappa_5, \kappa_6, \kappa_7, \kappa_8, \kappa_9,\kappa_{10} \right)
		\quad = \quad \\
	& \quad 
	(1, 1, 1, 1, 41774.858, 1, 2.5081*10^{-4}, 7.3335*10^{-3}, 1.1614*10^{-4}, 7.5610*10^{-5})~,
	\end{align*}
	\end{small}
there are two steady states:
	\begin{align*}
	x^* \quad &= \quad ( 63.143335, 136.35902, 41577.356) \quad \quad {\rm and } \\	
	x^{**} \quad &= \quad ( 25473.839, 1007.5644, 15295.454)~.
	\end{align*}
\end{example}

% SUBSECTION:
% CITE BADAL's ONE-REACTION RESULT
\subsection{Classification of multistationary one-reaction CFSTRs} \label{sec:1Rxn}
We now recall the following theorem, due to the first author:
\begin{theorem}[\cite{1rxn}] \label{thm:1rxn}
\begin{enumerate}
	\item Consider a CFSTR which contains only one non-flow reaction:
	\begin{align*}
	 a_1 X_1 + a_2 X_2 + \cdots + a_s X_s \quad \ra \quad  b_1 X_1 + b_2 X_2 + \cdots + b_s X_s~,
	\end{align*}
	where $a_i,b_i\geq 0$.  
	Then the CFSTR admits multiple positive \ma steady states if and only if $\sum_{i:~b_i>a_i} a_i >1$.
	Moreover, these multistationary CFSTRs admit nondegenerate steady states.
	\item Consider a CFSTR in which the only non-flow reactions consist of a pair of reversible reactions:
	\begin{align*}
	 a_1 X_1 + a_2 X_2 + \cdots + a_s X_s \quad \lra \quad  b_1 X_1 + b_2 X_2 + \cdots + b_s X_s~,
	\end{align*}
	where $a_i,b_i \geq 0$.  
	The CFSTR admits multiple positive \ma steady states if and only if the following holds:
	\begin{align*}
	\sum_{i:~b_i>a_i} a_i >1 \quad {\rm or} \quad \sum_{i:~a_i>b_i} b_i >1~.
	\end{align*}
	Moreover, these multistationary CFSTRs admit nondegenerate steady states.
\end{enumerate}
%Informally:
% A one-reaction network admits multiple steady states if and only if it has one of the following `subreactions': $2A\rightarrow3A$ or $A+B\rightarrow2A+2B$.
\end{theorem}
The current work was motivated by the question of whether a similar theorem exists for the class of CFSTRs that consists of networks with {\em two} reversible nonflow reactions and their sub-CFSTRs.

% SUBSECTION: how JC used (in `Simplifying' paper)
%\subsection{Ruling out multistationarity by the Jacobian Criterion}\label{sec:JC}
%We now review the content of recent work of the authors~\cite{Simplifying}.  [To do: explain the JC, and explain how recent work simplifies it.]  We need to recall the following definition:
%DEFINITION: Total Molecularity of a species
%\begin{definition}
%Consider a \CRN whose set of non-flow reactions is given by:
%	\begin{align*}
%	y_1 & \ra y_1' \quad & y_2 & \ra y_2' \quad & \ldots & \quad & y_l & \ra y_l'  \\
%	y_{l+1} & \lra y_{l+1}' \quad & y_{l+2} & \lra y_{l+2}' \quad & \ldots & \quad & y_{l+k} & \lra y_{l+k}' ~,
%	\end{align*}
%where none of the first $l$ reactions is reversible.  The {\bf total molecularity} of a species $X_i$ in the network, denoted by $\op{TM}(X_i)$, is the following integer:
%	\begin{align*}
%	\op{TM}(X_i) \quad = \quad \sum_{j=1}^{l+k} \left( y_{ji} + y_{ji}' \right)~,
%	\end{align*}	
%	where $y_{ji}$ is the stoichiometric coefficient of species $i$ in the complex $y_j$.
%\end{definition}

% THM: TM at most 2
%\begin{theorem}
%\label{thm:TM}
%If all species of a network have total molecularity less than or equal to two, then the network passes the Jacobian Criterion.  If in %addition, the network is a CFSTR network or is weakly-reversible, then it does not admit multiple steady states.
%\end{theorem}

%----------------------------------
% SECTION: Lifting MSS from subnetworks
%----------------------------------
\section{Lifting multistationarity from subnetworks} \label{sec:lift}

Consider the following question: {\em if a subnetwork $N$ of a network $G$ admits multiple positive steady states, then does $G$ as well?}  Theorem~\ref{thm:subnetwork} asserts that the answer to this question is `yes,' provided that the steady states are nondegenerate and the two networks share the same stoichiometric subspace (note that the stoichiometric subspace of $N$ is always contained in that of $G$).  The proof lifts each steady state $x^*$ of $N$ to a nearby steady state of $G$.

% THEOREM: lifting MSS from subnetworks
\begin{theorem} \label{thm:subnetwork}
Let $N$ be a subnetwork of a \CRN $G$ such that they have the same stoichiometric subspace: $\Sbar_N = \Sbar_G$. Let $\famKins$ be a parametrized family of kinetics on the species of $G$. Then the following holds: 
	\begin{itemize}
	\item If $N$ admits multiple nondegenerate positive $\famKins$ steady states, then $G$ does as well.  Additionally, if $N$ admits finitely many such steady states, then $G$ admits at least as many.  
	\item Moreover, if $N$ admits multiple positive {\em exponentially stable} steady states, then $G$ does as well.  Additionally, if $N$ admits finitely many such steady states, then $G$ admits at least as many.
	\end{itemize}
\end{theorem}

We note that our theorem is similar to Theorem~2 in work of Craciun and Feinberg~\cite{ME_entrapped}; their theorem % is stated for more general kinetics than mass-action, and 
allows multiple steady states to be lifted from an `entrapped species' network (that is, only certain species are in the outflow) to the corresponding `fully diffusive' network (all species are in the outflow).  In addition, their theorem is stated as a contrapositive version of ours. % The proof of Theorem~\ref{thm:subnetwork} is a straightforward extension of the proof of Craciun and Feinberg.
% As an alternative to using the Implicit Function Theorem as was done by Craciun and Feinberg, o
Our proof of Theorem~\ref{thm:subnetwork} makes use of the following homotopy theory result, which is a modified form of Theorem~1.1 in Craciun, Helton, and Williams \cite{CHW08}.
% Homotopy lemma
	\begin{lemma} \label{lem:homotopy}
	Let $S \subset \R^n$ be a vector subspace, let $\invtPoly \subset \R^n$ be a polyhedron contained in an affine translation of $S$, and let $\Omega \subset \op{int}(\invtPoly)$ be a bounded domain in the relative interior of $\invtPoly$. Assume that $g_{\lambda}: \overline{\Omega} \to S$, for $\lambda \in [0,1]$, is a continuously-varying family of smooth functions such that
	\begin{enumerate}
	\item for all $\lambda \in [0,1]$, $g_{\lambda}$ has no zeroes on the boundary of $\Omega$, and 
	\item for $\lambda=0$ and $\lambda=1$,  $\im dg_{\lambda} (\overline x) = S $ for all $x \in \Omega$.
	\end{enumerate}	 
Then the number of zeroes of $g_0$ in $\Omega$ equals the number of zeroes of $g_1$ in $\Omega$.	
	\end{lemma}

% PROOF of Theorem for lifting MSS from subnetworks
We now prove Theorem~\ref{thm:subnetwork}.
\begin{proof}[Proof of Theorem~\ref{thm:subnetwork}]
First, note that the network $G$ and its subnetwork $N$ must have the same set of species $\SSS$ in order for their stoichiometric subspaces to coincide.  We let $\Sbar$ denote the shared stoichiometric subspace: $\Sbar := \Sbar_G = \Sbar_N $.  Now, let $\RR'$ denote the set of reactions of $G$ that are not in $N$:
$	\RR_G \quad = \quad \RR_N ~\dot \cup~ \RR' $.
We now assume that the subnetwork $N$ admits multiple nondegenerate positive steady states; that is, there exist rate constants $\kappa_1^*, \kappa_2^*, \dots, \kappa_{|\RR_N|}^* \in \Rplus$ such that there exist distinct, stoichiometrically compatible, nondegenerate positive steady states $x^*$ and $x^{**}$ of the \CRS $(N,\kappa_{i}^*)$ arising from $\famKins$.  Write $f_N$ for the differential equations of $(N,\kappa_{i}^*)$. Now $x^*$ is a nondegenerate steady state, % i.e. f_N(x^*)=0 and $\im d_x f_N (x^*)=\Sbar$
 so there exists a relatively open ball $\Omega$ around $x^*$ in the interior of $\invtPoly$ such that (1) $x^*$ is the unique steady state (zero of $f_N$) in $\Omega$, and (2) $\im d f_N (x)=\Sbar$ for all $x \in \Omega$.  Note that (2) can be accomplished because the non-vanishing of a determinant is an open condition and because the matrix $df_N(x)$ varies continuously in $x$.
 
 For any vector of reaction parameters $\kappa \in \Rplus^{|\RR'|}$, we define the following the following family of functions for $0 \leq \lambda \leq 1$:
	\[
	g^{\kappa} _{\lambda} (x)  := f_N(x) + \sum_{k \in \RR'} (y_k'-y_k) R^{\lambda \kappa_k}_{k} (x)~.
	\] 
 It follows that $g^{\kappa} _{\lambda}(x)$ gives the differential equations~\eqref{eq:main_general} of the \CRS arising from the network $G$ and the following reaction parameters with respect to the kinetics $\famKins$:
	\begin{align} \label{eq:G_rates}
	(\kappa^*,\kappa) \quad := \quad \left( \kappa_1^*, \kappa_2^*, \dots, \kappa_{|\RR_N|}^*, \lambda \kappa_1, \lambda \kappa_2, \dots, \lambda \kappa_{|\RR'|} \right)\in \Rplus^{|\RR_G|}~.
	\end{align} 
Note that $g^{\kappa} _{0} (x) = f_N(x)$.  Next, by continuity in $\kappa$ and the compactness of the boundary of $\Omega$, there exists a vector of reaction parameters $\kappa^{\dagger} \in \Rplus^{|\RR'|}$ such that for all $0 \leq \lambda$, the function $g^{\kappa^{\dagger}} _{\lambda} (x)$ has no zeroes on the boundary of $\Omega$.  By continuity in $\lambda$, and by scaling $\kappa^{\dagger}$ smaller if necessary, we may assume additionally that $\im d_x g^{\kappa^{\dagger}} _{\lambda} (x)=\Sbar$ for all $x \in \Omega$.  Therefore, Lemma~\ref{lem:homotopy} allows us to conclude that the \CRS $\left(G,(\kappa^*,\lambda\kappa^{\dagger}) \right)$ has a nondegenerate steady state in the ball $\Omega$ for all $0 \leq \lambda \leq 1$.

We now complete the proof by repeating the argument with $\xTwo$, taking care that the ball around $\xTwo$ does not intersect that of $\xOne$; we replace $\kappa^{\dagger}$ by a scaled-down version ($\mu \kappa^{\dagger}$ for some $0 < \mu < 1$) if necessary.  It follows that $\left(G,(\kappa^*,\lambda\kappa^{\dagger}) \right)$ has at least two nondegenerate steady states.  The case of three or more nondegenerate steady states generalizes in a straightforward way. 
%STABILITY
For the stability result, we simply note that the eigenvalues of a matrix vary continuously under continuous perturbations (in this case, arising from the parameter $\lambda$). 
\end{proof}

% REMARK: justify Siegal-Gaskins et al.
\begin{remark} \label{rmk:justify}
One application of Theorem~\ref{thm:subnetwork} is that it provides the mathematical justification for the technique of Siegal-Gaskins {\em et al.}\ which establishes bistability in the mass-action setting by way of `network ancestry'~\cite{Emergence}.  In their examination of $40,680$ small gene regulatory networks, $14,721$ initially were established to be bistable by the implementation of Advanced Deficiency Theory in the CRN Toolbox~\cite{Toolbox}, and an additional $22,050$ were classified as bistable by virtue of containing one of the $14,721$ bistable networks as a subnetwork (`ancestor') such that both networks have the same stoichiometric subspace.  A similar approach is taken by Conradi {\em et al.}\ for lifting multiple steady states from certain subnetworks called `elementary flux modes'~\cite{Conradi_subnetwork}.  We note that their criterion for lifting steady states does not require that the stoichiometric subspaces of the network and its subnetwork to coincide \cite[Supporting Information]{Conradi_subnetwork}. %, whereas Theorem~\ref{thm:subnetwork} here does not.
\end{remark}

% EXAMPLE: why the hypotheses in the Theorem are necessary
The next example illustrates why the hypothesis of nondegeneracy %and the coincidence of the stoichiometric compatibility classes are
is required in Theorem~\ref{thm:subnetwork}.  A larger such example appears in the work of Craciun and Feinberg \cite[\S 6]{ME_entrapped}.
\begin{example}
Consider the following (non-CFSTR) network:
% TRIANGLE NETWORK PLUS ONE REVERSIBLE REACTION
	\begin{align} \label{eq:deg_net}
	\begin{xy}<10mm,0cm>:
	%TRIANGLE points
	(0,1.3)                  ="A+B"  *+!D{ B }  *{}; % TOP
	(-2,0)                  ="D"  *+!R{A}  *{}; % LEFT
	(2,0)                  ="C"  *+!L{C}  *{}; % RIGHT
	% One more reversible reaction: 2 complexes
	(4,0.5) = "A+C" *+!R{ A+C }  *{}; % Left complex
	(6,0.5) = "2B" *+L{~~~~~2B} *{}; % Right complex
	%Rate constants (none)
	%Reaction arrows
	{\ar "A+B"+(-0.15,0)*{};"C"+(-0.3,.15)*{}  };
	{\ar "C"+(0,0.15)*{};"A+B"+(0.15,0)*{}};      
	{\ar "A+B"+(-0.45,-0.15)*{};"D"+(-0.15,0.15)*{}};
	{\ar "D"+(0.15,0.15)*{};"A+B"-(0.15,0.15)*{}};
	{\ar "D"+(0.15,-0.10)*{};"C"+(-0.15,-0.10)*{}};     
	{\ar "C"+(-0.15,0.05)*{};"D"+(+0.15,0.05)*{}};  
   	%Reaction arrows for reversible reactions
	{\ar "A+C"+(0.15,-0.10)*{};"2B"+(-0.15,-0.10)*{}};     
	{\ar "2B"+(-0.15,0.05)*{};"A+C"+(+0.15,0.05)*{}};  		
	\end{xy}
	\end{align}
The CRN Toolbox~\cite{Toolbox} determines that network~\eqref{eq:deg_net} does not admit multiple positive mass-action steady states, but the following subnetwork does admit multiple degenerate positive steady states:
	\begin{align} \label{eq:deg_sub}
	A  \stackrel[~]{\kappa_{1}}{\leftarrow} B  \stackrel[~]{\kappa_{2}}{\ra} C \quad \quad  A+C  \stackrel[~]{\kappa_{3}}{\ra} 2B~.
	\end{align}
In fact, it is straightforward to verify that steady states exist for network~\eqref{eq:deg_sub} if and only if $\kappa_1=\kappa_2$, and in this case, each two-dimensional compatibility class contains an infinite one-dimensional set of degenerate steady states.
\end{example}

% ALSO: note that the converse of main theorem does not hold (ex: G has MSS, N does not)

One way for a network and its subnetwork to share the same stoichiometric subspace is for the subnetwork to be obtained by making some reversible reactions irreversible.  Thus, Theorem~\ref{thm:subnetwork} yields the following corollary.  
% COROLLARY: lift MSS in the case of making irreversible reactions reversible
\begin{corollary} \label{cor:lift_irr}
For a \CRN $N$, let $\famKins$ be a parametrized family of kinetics on the species of $N$.  Let $G$ be a network obtained from $N$ by making some irreversible reactions of $N$ reversible.  Then if $N$ admits multiple nondegenerate positive $\famKins$ steady states, then $G$ does as well.
\end{corollary}

The next corollary states that Theorem~\ref{thm:subnetwork} allows multiple positive steady states to be lifted from a sub-CFSTR to a \fd CFSTR.  % without requiring that the stoichiometric subspaces of the two networks are the same.  
Therefore, the set of minimal multistationary CFSTRs (with respect to the subnetwork relation) completely defines the set of all \fd multistationary CFSTRs: {\em a \fd CFSTR admits multiple steady states if and only if it contains as a subnetwork one of these minimal CFSTRs}.  This result will be useful in our classification of small multistationary CFSTRs in Section~\ref{sec:class}.
% COROLLARY: lift MSS in the case of CFSTR
\begin{corollary} \label{cor:liftCFSTR}
Let $N$ be a sub-CFSTR of a \fd CFSTR $G$, and let $\famKins$ be a parametrized family of kinetics on the species of $G$.  Then, if $N$ admits multiple nondegenerate positive $\famKins$ steady states, then $G$ does as well.
\end{corollary}
\begin{proof}
Assume that the species of $N$ are $X_1, X_2, \dots, X_{s_1}$ and the species of $G$ are $X_1, X_2, \dots,$ $X_{s_1+s_2}$.  Let $N'$ be the CFSTR obtained from $N$ by appending the flow reactions $0 \lra X_{s_1+1}$, $0 \lra X_{s_1+2}$, \dots, $0 \lra X_{s_1+s_2}$ for all species of $G$ that are not in $N$.  Clearly, $N'$ is a subnetwork of $G$, and they share the same stoichiometric subspace, namely, $\R^{s_1+s_2}$.  By applying Theorem~\ref{thm:subnetwork} to $N'$ and $G$, we see that if $N'$ admits multiple nondegenerate positive steady states, then $G$ does as well.  Therefore, it remains only to show that $N$ admits multiple nondegenerate positive steady states if and only if $N'$ does.  

Consider any outflow rate parameter $\kappa_{s_1+i}^{\rm out}>0$ for one of the new outflow reactions.  Then by Assumption~\ref{assump:rates}, the rate function $R_{\kappa_{s_1+i}^{\rm out}}^{X_{s_1+i}\to 0}(x)$ depends only on $x_{s_1+i}$ and is increasing in $x_{s_1+i}$ from 
	\[
	0=R_{\kappa_{s_1+i}^{\rm out}}^{X_{s_1+i}\to 0}(x_1,\dots, x_{s_1+i-1},0,x_{s_1+i+1}, \dots x_{s_1+s_2})~.
	\]
 As for the corresponding inflow rate function, Assumption~\ref{assump:rates} implies that $R_{\kappa_{s_1+i}^{\rm in}}^{0 \to X_{s_1+i}}(x)$ is a positive constant function, and this constant depends only on the parameter ${\kappa_{s_1+i}^{\rm in}}$ and is in fact increasing in this parameter for sufficiently small values, with $R_{0}^{0 \to X_{s_1+i}}(x)=0$.  Thus, we can choose a sufficiently small inflow parameter ${\kappa_{s_1+i}^{\rm in}} > 0$ such that there exists a positive value $x^*_{s_1+i}>0$ for which the rate functions are equal at $x$ when $x_{s_1+i}=x^*_{s_1+i}$.

Therefore, it follows that 
	% Steady state of N
	$
	x^*  =  \left( x^*_1, x^*_2, \dots, x^*_{s_1} \right) \in \Rplus^{s_1}
	$
is a nondegenerate positive steady state of the system $\left(N,(\kappa_1, \kappa_2, \dots, \kappa_{|\RR_N|})\right)$ if and only if the concentration vector 
	% Steady state of N'
	\begin{align*}
	\left( x^*_1, x^*_2, \dots, x^*_{s_1}, 
	x^*_{s_1+1}, \dots, x^*_{s_1+s_2}
%		\frac{\kappa_{s_1+1}^{\rm in} }{ \kappa_{s_1+1}^{\rm out}}, 
%		\frac{\kappa_{s_1+2}^{\rm in} }{ \kappa_{s_1+2}^{\rm out}},  \dots , 
%		\frac{\kappa_{s_1+s_2}^{\rm in} }{ \kappa_{s_1+s_2}^{\rm out}}  
	\right) \in \Rplus^{s_1+s_2}
	\end{align*}
is a nondegenerate positive steady state of the system 
		\[
	\left(N',(\kappa_1, \kappa_2, \dots, \kappa_{|\RR_N|}, \kappa_{s_1+1}^{\rm in},  \kappa_{s_1+1}^{\rm out} \dots,  \kappa_{s_1+s_2}^{\rm out} ) \right)~,
	\]
	 where the rates $\kappa_{s_1+i}^{\rm in}$ and $\kappa_{s_1+i}^{\rm out}$ are chosen as described above.    This completes the proof.
\end{proof}

%----------------------------------
% SECTION: Lifting MSS from embedded networks (in case of CFSTR)
%----------------------------------
\section{Lifting mass-action multistationarity from embedded CFSTRs} \label{sec:lift_emb}
Corollary~\ref{cor:liftCFSTR} stated that multistationarity can be lifted from sub-CFSTRs; in this section, we generalize the result to the case of embedded CFSTRs  in the mass-action setting (Corollary~\ref{cor:liftEmbCFSTR}).

We first need the following generalization of inflow/outflow reactions in order to allow for reactions such as $A \lra 2A$ which also have a mass-action steady state at $\overline{x_A}=1$ when the two reaction rate constants are equal.
% FLOW-TYPE SUBNETWORK
\begin{definition} \label{def:flow-type}
A {\em mass-action flow-type subnetwork} for a species $X_i$ of a \CRN $G$ is a nonempty subnetwork $N$ of $G$ such that 
	\begin{enumerate}
	\item the reactions in $N$ involve only species $X_i$, and 
	\item there exists a choice of reaction rate constants $\kappa_r^*$ for the reactions $r \in \RR_N$ of $N$ such that for the resulting mass-action system of this subnetwork $N$, $\overline{x_i}= 1$ is a nondegenerate steady state.

	\end{enumerate}
\end{definition}

The following theorem is analogous to Theorem~\ref{thm:subnetwork}.
% THEOREM: lifting MSS from embedded networks
\begin{theorem}  \label{thm:emb}
Let $N$ be an embedded network of a network $G$ such that
	\begin{enumerate}
	\item the stoichiometric subspace of $N$ is full-dimensional: $\Sbar_N = \R^{|\SSS_N|}$, and 
	\item for each species $X_i$ that is in $G$ but not in $N$, there exists a mass-action flow-type subnetwork of $G$ for $X_i$. 	 %makes G full-dim'l
	\end{enumerate}
Then the following holds:
	\begin{itemize}
	\item If $N$ admits multiple nondegenerate positive mass-action steady states, then $G$ does as well.  Additionally, if $N$ admits finitely many such steady states, then $G$ admits at least as many.  
	\item Moreover, if $N$ admits multiple positive {\em exponentially stable} mass-action steady states, then $G$ does as well.  Additionally, if $N$ admits finitely many such steady states, then $G$ admits at least as many.
	\end{itemize}
\end{theorem}

The proof of Theorem~\ref{thm:emb} requires the following lemma, which states that for certain simple embedded networks obtained by removing only one species, each nondegenerate steady state $u$ can be lifted to a steady state of the larger network that is near $(u,1)$.
% LEMMA: just remove one species, and there are no additional reactions except flow-type
\begin{lemma} \label{lem:emb}
Let $G$ be a \CRN with $s$ species denoted by $X_1, X_2, \dots, X_s$, and let $N$ be an embedded network of $G$ with $s-1$ species $X_1, X_2, \dots, X_{s-1}$ such that $N$ is full-dimensional: $\Sbar_N = \R^{s-1}$.  Assume that the reactions of $G$ and the reactions of $N$ can be written as, respectively, $\RR_G=\{ \widetilde{R_1}, \widetilde{R_2}, \dots, \widetilde{R_m}, R_{m+1},  \dots, R_{m+n} \} $ and $\RR_N=\{ R_1, R_2, \dots, R_{m} \}$ such that:
	\begin{enumerate}
	\item for $i=1,2,\dots, m$, the reaction $R_i$ of $N$ is obtained from the corresponding reaction $\widetilde{R_i}$ of $G$ by removing species $X_s$, and
	\item  all remaining reactions of $G$, namely $\{R_{m+1},R_{m+2}, \dots, R_{m+n} \}$, together form a mass-action flow-type subnetwork for the species $X_s$.
	\end{enumerate}
For a choice of rate constants $\kappa^*_1, \kappa^*_2, \dots, \kappa^*_m>0$, let 
$\Sigma\left(N,\{\kappa^*_1, \kappa^*_2, \dots, \kappa^*_m\}\right)$ 
denote a finite set of nondegenerate positive mass-action steady states of the system arising from $N$ and the $\kappa^*_i$.  Then for sufficiently small $\ep > 0$, there exist reaction rate constants $\kappa^*_{m+1}, \kappa^*_{m+2}, \dots, \kappa^*_{m+n}$ for the flow-type subnetwork of $G$ such that for all $u \in \Sigma\left(N,\{\kappa^*_1, \kappa^*_2, \dots, \kappa^*_m\}\right)$, there exists a nondegenerate positive mass-action steady state $\widetilde{u}$ % \in \Sigma\left(G,\{\kappa^*_1, \kappa^*_2, \dots, \kappa^*_{m+n}\}\right)$
 of the system arising from $G$ and $\kappa^*_1, \kappa^*_2, \dots, \kappa^*_{m+n}$ with $| \widetilde{u} - (u,1) | < \ep$.  Additionally, if $u$ is exponentially stable, then $\widetilde{u}$ is as well.
\end{lemma}
% PROOF
\begin{proof}
Fix a choice of reaction rate constants $\kappa^*_1, \kappa^*_2, \dots, \kappa^*_m$, and let $\SigSet:=\Sigma\left(N,\{\kappa^*_1, \kappa^*_2, \dots, \kappa^*_m\}\right)$ be as in the statement of the lemma.

% G - disjoint union of two subnetworks
We view $G = \widetilde{N} \dot{\cup} M$ as the disjoint union of two subnetworks, one which consists of the reactions $\widetilde{R_1}, \widetilde{R_2}, \dots, \widetilde{R_m}$, which we denote by $\widetilde{N}$, and the second which consists of $R_{m+1}, R_{m+2}, \dots,$ $R_{m+n}$, which we denote by $M$.  As $M$ is a mass-action flow-type subnetwork for the species $X_s$, there exist rate constants $\overline{\kappa_{m+1}}, \overline{\kappa_{m+2}}, \dots, \overline{\kappa_{m+n}} > 0$ such that the resulting mass-action ODE system, denoted by $f_M(x_s)$, has a nondegenerate steady state at $\overline{x_s} =1$.

% write down the ODE systems
Next, we denote by $f_{\widetilde{N}} (x)$ the mass-action ODE system~\eqref{eq:main} arising from the subnetwork $\widetilde{N}$ and the fixed rate constants  $\kappa^*_1, \kappa^*_2, \dots, \kappa^*_m$.  
Consider the following map from $\Rnn \times \Rnn^{s}$ to $\R^s$:
	\begin{align} \label{eq:G-sys}
	f_G(k,x)  :=  
		\left(
			f_{\widetilde{N},1} (x), ~ f_{\widetilde{N},2} (x), ~ \dots~, f_{\widetilde{N},s-1} (x), ~
				f_{\widetilde{N},s} (x)  +  k f_M(x_s)  
		\right)~.
	\end{align}
(Note that 	$f_{\widetilde{N},i}$ denotes the $i$-th coordinate function of $f_{\widetilde{N}}$.)  It follows that $f_G(k,x)$ denotes the mass-action ODEs for the network $G$ with respect to the rate constants $$\kappa_1^*, \kappa_2^*, \dots , \kappa_m^*, k \overline{\kappa_{m+1}}, k  \overline{\kappa_{m+2}}, \dots k \overline{\kappa_{m+n}}~.$$  
% Rescale our last coordinate
We scale the last coordinate of $f_G(k,x)  $ by $1/{k}$ and make the substitution $\delta = 1/k$ to obtain:
	\begin{align} \label{eq:G-sys-scaled}
	F_G(\delta,x)  :&=  
		\left(
			f_{\widetilde{N},1} (x), ~ \dots~,  f_{\widetilde{N},s-1} (x), ~
				  f_M(x_s)  
		\right) + 
		 \left( 0,\dots,0,  \delta f_{\widetilde{N},s} (x)   \right) \\
		&=: h(x) + \delta \left( 0,\dots,0,  f_{\widetilde{N},s} (x)   \right)~,  \notag
	\end{align}
where $h(x)$ is defined by the second equality.  	
Hence, it suffices to prove that for sufficiently small $\ep > 0$ and for all $u \in \SigSet$, there exists a $\delta>0$ such that there exists a nondegenerate zero $\widetilde{u}$ of $F_G(\delta,x)$ with $| \widetilde{u} - (u,1) | < \ep$.  

% h has nondeg zero a (u,1)
Fix $u \in \SigSet$.  We now claim that $h$ has a nondegenerate zero at $(u,1)$.  The final coordinate of $h$ satisfies $h_s(u,1)=f_M(1)=0$ by construction.  As for the remaining coordinates $i=1,2,\dots, s-1$, we compute
	\begin{align} \label{eq:h}
	h_i(u,1) = f_{\widetilde{N},i}(u,1) = f_N (u) = 0~.
	\end{align}
We now explain the second equality in~\eqref{eq:h}.  When the reaction $\widetilde{R_j}$ in $\widetilde{N}$ is $\widetilde{y_j} \ra \widetilde{y_j'}$ then the reaction $R_j$, given by $y_j \to y_j'$, is such that the projection of $\widetilde{y_j}$ onto the first $s-1$ coordinates is $y_j$ and similarly for $\widetilde{y_j'}$.  Thus, the reaction vector $\widetilde{y_j'} - \widetilde{y_j}$ projects to $y_j'-y_j$ and $(u_1,u_2,\dots, u_{s-1},1)^{\widetilde{y_j}} = (u_1,u_2,\dots, u_{s-1})^{y_j} $.  Finally, $(u,1)$ is nondegenerate, because $dh(u,1)$ is an $s \times s$-matrix in which the upper-left $(s-1) \times (s-1)$-submatrix  is the nonsingular matrix $df_N(u)$ and the bottom row is $(0,0,\dots,0, \frac{df_M}{dx_s}(1))$ with $ \frac{df_M}{dx_s}(1) \neq 0$ by hypothesis.

As $(u,1)$ is nondegenerate, there exists a constant $\ep(u)>0$ such that the resulting $\ep(u)$-neighborhood of $(u,1)$, which we denote by $\Omega$, is such that (1)~$\Omega$ is in the positive orthant $\Rplus^s$, (2) $(u,1)$ is the unique zero of $h$ in $\Omega$, and (3) $dh(x)$ is nonsingular for all $x \in \Omega$.  Consider again the function
	$F_G(\delta,x)$
defined in~\eqref{eq:G-sys-scaled}, and note that $F(0,x)=h(x)$.  By continuity in $\delta$ and the compactness of the boundary of $\Omega$, there exists $\delta(u)>0$ such that for all $0 \leq \delta \leq \delta(u)$, the function 
	$F_G(\delta,x)$ 
has no zeroes on the boundary of $\Omega$.  Again by continuity and by decreasing $\delta(u)$ if necessary, we may assume that $dF(\delta,x)$ (the matrix of partial derivatives with respect to the $x_1,x_2,\dots, x_s$) is nonsingular for all $0 \leq \delta \leq \delta(u)$ and for all $x \in \Omega$.  Therefore, Lemma~\ref{lem:homotopy} allows us to conclude that
	$F_G(\delta,x)$ 
has a unique nondegenerate zero $\widetilde{u}$ in $\Omega$ (that is, $| \widetilde{u} - (u,1) | < \ep(u)$) for all $0 \leq \delta \leq \delta(u)$.

Now let $\ep^*$ be the minimum of all such $\ep(u)$, where $u \in \SigSet$.  Additionally, we decrease $\ep^*$ if necessary so that the resulting $\ep$-neighborhoods of the points $u$ do not intersect.  The lemma now follows with the $\ep^*$ as a cut-off: given any $0 < \ep < \ep^*$, the above arguments for each $u$ can be made using $\ep$ in place of $\ep(u)$.  Taking the minimum, denoted by $\delta^*$, of the resulting cut-offs $\delta(u)$, we obtain nondegenerate zero $\widetilde{u}$ of $F_G(\delta^*,x)$ such that $| \widetilde{u} - (u,1) | < \ep$.  

For the stability result, the eigenvalues of a matrix vary continuously under continuous perturbations (in this case, arising from the parameter $\delta$). 
\end{proof}

% REMARK: time scales
\begin{remark}  \label{rmk:time-scale}
In the proof of Lemma~\ref{lem:emb}, the $s$-dimensional dynamical system~\eqref{eq:G-sys} may be represented by
\begin{align*}
\dot x = f_{\widetilde N} (x) + k g_M (x_s)~,
\end{align*}
where $g_M (x_s) := (0,\ldots,0,f_M(x_s))$. When $k$ is sufficiently large and $x_s$ is in the domain of attraction of $\overline {x_s}=1$ with $\abs{g_M(x_s)} \sim O(1)$, then $\dot x \sim kg_M(x_s)$, which has dynamics close to the one-dimensional system $\dot x_s = k f_M(x_s)$. However, when $x_s$ is close to $\overline {x_s}=1$ with $\abs{g_M(x_s)} \sim O(1/k^2)$, then $\dot x \sim f_{\widetilde N} (x)$, the dynamics of which are effectively those of an $(s-1)$-dimensional system. Thus by choosing $k$ large enough, we achieve a time-scale separation: on the fast time-scale, the dynamics are close to a one-dimensional system and on the slow time-scale, the dynamics are close to an $(s-1)$-dimensional system. Thus, we can lift the steady states from the smaller system to the full system. 
\end{remark}

We can now prove Theorem~\ref{thm:emb}.
% Proof of embedded theorem
\begin{proof}[Proof of Theorem~\ref{thm:emb}]
We begin by reducing to the case that $G$ has only one species that $N$ does not have: if $G$ has more than one additional species, we can lift multistationarity `one species at a time.'  Now denote the species of $N$ by $X_1, X_2, \dots, X_{s-1}$ and the species of $G$ by $X_1, X_2, \dots, X_s$.  Denote the reactions of $N$ by $y_1 \ra y_1', y_2 \ra y_2', \dots,  y_m \ra y_m'$, where $y_i, y_i' \in \mathbb{Z}_{\geq 0} ^ {s-1}$. As $N$ is an embedded network of $G$, we can write the reactions of $G$ as 
	$\RR_G=\{ \widetilde{R_1}, \widetilde{R_2}, \dots, \widetilde{R_m}, R_{m+1}, R_{m+2}, \dots, R_{m+n}, R_{m+n+1}, \dots, R_{m+n+p} \} $  such that:
	\begin{enumerate}
	\item for $i=1,2,\dots, m$, the reaction $R_i$ of $N$ is obtained from the corresponding reaction $\widetilde{R_i}$ of $G$ by removing species $X_s$, and
	\item  the reactions in $\{R_{m+1},R_{m+2}, \dots, R_{m+n} \}$ form a mass-action flow-type subnetwork for the species $X_s$.
	\end{enumerate}
%Define G'
We now let $G'$ denote the subnetwork of $G$ that consists of the reactions:  $ \widetilde{R_1}, \widetilde{R_2}, \dots, \widetilde{R_m},$ $R_{m+1}, R_{m+2}, \dots, R_{m+n} $.  Lemma~\ref{lem:emb} applies to this network $G'$ and its embedded network $N$, so $G'$ admits at least as many nondegenerate positive mass-action steady states as $N$ (and similarly for exponentially stable steady states).  Next, $G'$ is a subnetwork of $G$ that shares the same stoichiometric subspace (namely, $\R^s$), so by Theorem~\ref{thm:subnetwork}, $G$ admits at least as many nondegenerate positive mass-action steady states as $G'$ (and similarly for exponentially stable ones), so this completes the proof.
\end{proof}

% Why hypothesis 2 in theorem is necessary
We now illustrate the necessity of the hypothesis 2 of Theorem~\ref{thm:emb}. % that $G$ contains a flow-type subnetwork for the species not contained in the embedded network $N$. 

\begin{example} \label{ex:needHyp}
Consider the following (non-CFSTR) network $G$, which is adapted from a similar network that appears in work of Feinberg \cite{FeinbergMSSdefone}:
	\begin{align*}
	0 \stackrel[\kappa_2]{\kappa_1}{\rightleftarrows} A  \quad \quad 
	3A \stackrel[\kappa_4]{\kappa_3}{\rightleftarrows}  2A+B 
	\end{align*}
A straightforward calculation reveals that $G$ has a unique mass-action steady state, namely $(\overline{x_A},\overline{x_B})=\left(\frac{\kappa_1}{\kappa_2}, \frac{\kappa_1 \kappa_3}{\kappa_2 \kappa_4} \right)$. In fact, despite the fact that $A$ participates in a non-flow reaction, % (or a `true reaction'), 
the steady state value of $x_A$ is the same as it would be when considering only the flow subnetwork $0 \lra A$. Now consider the following embedded network $N$ obtained by removing the species $B$:
\begin{align*}
0 \lra A  \quad \quad
3A \lra 2A
\end{align*}
We see that $N$ satisfies the conditions of Theorem~\ref{thm:1rxn}, so $N$ admits multiple mass-action steady states.  Note that $N$ is an embedded network of $G$, but its multiple steady states can {\em not} be lifted to $G$; Theorem~\ref{thm:emb} does not apply because $G$ does not contain a flow-type subnetwork for the species $B$.

On the other hand, $N$ is an embedded network of the following network $G'$:
\begin{align*}
 B \lra  0 \lra A \quad \quad 
 3A \lra 2A +B
\end{align*}
which does contain  a flow-type subnetwork for the species $B$.  So, Theorem~\ref{thm:emb} does apply and thus we conclude that $G'$ admits multiple steady states.% by virtue of possessing the embedded multistationary network $M$:
\end{example}

We now have an analogue of Corollary~\ref{cor:liftCFSTR}.  
% COROLLARY: lift MSS in the case of CFSTR (for embedded networks)
\begin{corollary} \label{cor:liftEmbCFSTR}
Let $N$ be an embedded CFSTR of a \fd CFSTR $G$. %, and let $\famKins$ be a parametrized family of kinetics on the species of $G$ that is closed under scaling.  
Then, if $N$ admits multiple nondegenerate positive mass-action steady states, then $G$ does as well.
\end{corollary}

\begin{proof}
This follow directly from Theorem~\ref{thm:emb}, after noting that hypothesis 2 of the theorem is satisfied by the inflow/outflow reactions $0 \lra X_i$.
% was verified for \fd CFSTRs in the proof of Corollary~\ref{cor:liftCFSTR} for the pair of reversible inflow/outflow reactions $0 \lra X_{i}$.
\end{proof}

%------------------
% SECTION: Atoms
%------------------
\section{CFSTR atoms of multistationarity} \label{sec:atom}
In the previous section, we saw that a CFSTR is multistationary in the mass-action setting if and only if an embedded CFSTR is multistationary; now we call the minimal such networks  `atoms of multistationarity.'  In Section~\ref{sec:class}, we will classify certain two-reaction atoms of multistationarity (see Corollary~\ref{cor:class}).  % that the two-reaction bimolecular CFSTRs are determined completely by the 11 minimal (with respect to the `embedded' relation) multistationary CFSTRs.  
% Note that in the remainder of this work, we restrict our attention to mass-action kinetics.

% ATOMS - definition
\begin{definition} \label{def:atom}
	\begin{enumerate}
	\item
	A \fd CFSTR is a {\em CFSTR atom of multistationarity} if it admits multiple nondegenerate positive mass-action steady states and it is minimal with respect to the embedded network relation among all such \fd CFSTRs.
	%\fd CFSTRs that admit multiple nondegenerate positive mass-action steady states.
	\item 
	A \fd CFSTR $G$ is said to {\em possess a CFSTR atom of multistationarity} if there exists an embedded network $N$ of $G$ that is a CFSTR atom.
	\end{enumerate}
\end{definition}
We now restate Corollary~\ref{cor:liftEmbCFSTR} in the following way, which motivates the above definition and suggests that compiling a list of atoms is desirable.
% FORMER CONJECTURE
	\begin{corollary} \label{cor:atom}
	A \fd CFSTR possesses a CFSTR atom of multistationarity if and only if it admits multiple nondegenerate positive mass-action steady states.
	\end{corollary}
\begin{proof}
The reverse direction is clear: a multistationary CFSTR is either itself a CFSTR atom of multistationarity or contains one.  The forward direction is Corollary~\ref{cor:liftEmbCFSTR}.
\end{proof}	
We also can rephrase Theorem~\ref{thm:1rxn} in the following way:
% and part 2 of Corollary~\ref{cor:class} in the following way:
% COROLLARY
	\begin{corollary} \label{cor:atom2}
	A one-reaction CFSTR is a CFSTR atom of multistationarity if and only if it consists of one non-flow reaction and that non-flow reaction has one of the following two forms:
	\begin{align} \label{eqn:1rxnatom}
	a_1 X \quad \ra \quad a_2 X~, \quad {\rm or} \quad
	X+ Y \quad \ra \quad b_1 X + b_2 Y~, 
	\end{align}
	where $a_2>a_1>1$, or, respectively, $b_1>1$ and $b_2 >1$.   A one-reaction CFSTR possesses one such CFSTR atom of multistationarity if and only if it admits multiple nondegenerate positive mass-action steady states.
%	Among bimolecular two-reaction CFSTRs, there are 11 CFSTR atoms of multistationarity (those displayed in bold/red in Figure~\ref{fig:35}).  A bimolecular two-reaction CFSTR possesses one of these 11 CFSTR atoms of multistationarity if and only if it admits multiple positive steady states.
	\end{corollary}

%{\bf Note}: This is a weaker concept of an atom than is the appropriate sense for one-reaction CFSTRs: $2A \rightarrow 3A$ and $3A \rightarrow 4A$ would be considered distinct atoms. {\it Can the characterization be made stronger so that there are fewer atoms?}

%I wonder if going through the list of atoms, but in fact listing real atoms this time, might help find a pattern? 
% What would the algorithm be to discover atoms in the general CRN (not-necessarily-CFSTR) setting?

We end this section by posing the following questions:
% QUESTIONS ABOUT ATOMS
\begin{enumerate}
\item {\em Is there a good characterization of CFSTR atoms of multistationarity?} For instance, even though there are countably infinitely many one-reaction CFSTR atoms, Corollary~\ref{cor:atom2} gives a simple characterization of all such one-reaction atoms. In particular, a one-reaction atom contains at most two species, and furthermore each of these atom types is characterized by exactly two parameters, $(a_1,a_2)$ or $(b_1,b_2)$ in equation \eqref{eqn:1rxnatom}. 
\item {\em Is there a good notion of `atom of multistationarity' outside of the CFSTR setting?}  If so, then a CFSTR atom might contain as an embedded network, a more general atom, which is obtained by removing some flow reactions and possibly more reactions.  For example, we can remove the outflow reaction $A \ra 0$ from the CFSTR atom arising from $A\ra 2A ~~ A+B \ra 0$ (see the top of Figure~\ref{fig:35} in the next section) and maintain multistationarity, but removing $B \ra 0$ destroys multistationarity.
\end{enumerate}
Beginning in the next section, we will give a partial answer to the first question above for two-reaction CFSTRs.

%----------------------------------
% SECTION: Enumeration
%----------------------------------
\section{Enumeration of reversible bimolecular two-reaction CFSTRs} \label{sec:enum}
% AIM: answer one question
The remainder of this work is dedicated to answering the following question:
\begin{question} \label{q:2RxnMSS}
Which bimolecular two-reaction \fd CFSTRs admit multiple positive mass-action steady states?
\end{question}
% Define 'bimolecular' and 'two-reaction'
By {\em bimolecular} we mean that each complex contains at most two molecules: the complexes $0$, $A$, $2A$, and $A+B$ are permitted, but $2A+B$ is not. A {\em two-reaction CFSTR} refers to a CFSTR in which the non-flow reactions consist of two pairs of reversible reactions, one reversible reaction and one irreversible reaction, or two irreversible reactions.  For instance, the three CFSTRs~\eqref{eq:running}, \eqref{eq:running_sub}, and \eqref{eq:running_emb} in Example~\ref{ex:running} are among the bimolecular two-reaction CFSTRs for which we would like to answer Question~\ref{q:2RxnMSS}.  Let us note that reactions of the form $0 \lra 2A$ or $0 \lra A+B$ (or any of the directed versions) are considered non-flow reactions. Finally, if we define two networks to be {\em equivalent} if there exists a relabeling of the species that transforms the first network into the second network, we aim to list only one network from each such equivalence class.  For example, the two CFSTRs in which the non-flow reactions are $2A \la A+B \la A$ and $C \ra B+C \ra 2C$, respectively, are both in the same equivalence class.

% Use Corollary
Note that it is sufficient to enumerate the possible non-flow subnetworks of our CFSTRs of interest; for example, if the non-flow subnetwork is 
	\begin{align} \label{eq:dir_sub}
	2A \la A+B \la A~,
	\end{align}
then the corresponding CFSTR is obtained by including the flow reactions for species $A$ and $B$.  In addition, Corollary~\ref{cor:lift_irr} implies that a non-reversible CFSTR (for example, the one arising from~\eqref{eq:dir_sub}) does not admit multiple nondegenerate positive steady states if the corresponding reversible CFSTR (for example, the one arising from $2A \lra A+B \lra A$) does not.  Therefore, we will proceed to answer Question~\ref{q:2RxnMSS} by completing the following steps:
	% 3 steps to complete
	\begin{enumerate}
	\item Enumerate all {\em reversible} bimolecular two-reaction networks.
	\item Determine which of the \fd CFSTRs arising from networks in Step~1 admit multiple  positive  mass-action steady states.
	\item Of those reversible CFSTRs that admit multiple positive steady states which were found in Step 2, determine which sub-CFSTRs admit multiple positive steady states. 
	\end{enumerate}
The current section describes how we performed Step 1 (see Algorithm~\ref{algor:enum}), and in Section~\ref{sec:class}, we explain how we completed Steps 2 and 3.

%------------------------
% SUBSECTION: Partitions
\subsection{The total molecularity partition of a \CRN}
We now explain how a network defines a `total molecularity partition'; two-reaction networks will be enumerated by these partitions in Algorithm~\ref{algor:enum}.  Recall that a {\em partition} of a positive integer $m$ is an unordered collection of positive integers that sum to $m$; by convention, we write the partition as $(m_1,m_2, \cdots, m_n )$, where the {\em parts} $m_i$ are weakly decreasing: $m_1 \geq m_2 \geq \cdots \geq m_n \geq 1$.  Partitions of $m=4,5,6,7,8$ are listed (partially) in Table~\ref{table:ptn}.  
%----------
%:Table 1
\begin{table}
\begin{center}
% Caption must be _above_ table (for journal style file)
\caption{Here we (partially) list the partitions of $m=4,5,6,7,8$ in lexicographic order.  The numbers of partitions, which are known as the Bell numbers, are listed in the last column.  \label{table:ptn} }
\begin{tabular}{ |c|c|c|}
\hline
$m$ 		& Partitions of $m$	& \# of partitions	 \\
\hline
4 &  (4), (3, 1), (2, 2), (2, 1, 1), (1, 1, 1, 1)
		& 5  \\
\hline
5 &  (5), (4, 1), (3, 2), (3, 1, 1), (2, 2, 1), (2, 1, 1, 1), (1, 1, 1, 1 , 1)
		& 7  \\
\hline
6 &  (6), (5, 1), (4, 2), (4, 1, 1), (3, 3), (3, 2, 1),~ % (3, 1, 1, 1), (2, 
  %2, 2), (2, 2, 1, 1), (2, 1, 1, 1, 1), 
	\dots~, (1, \dots , 1)
		& 11  \\
\hline
7 & (7), (6, 1), (5, 2), (5, 1, 1), (4, 3), (4, 2, 1), %(4, 1, 1, 1),(3, 
 % 3, 1), (3, 2, 2),
	 \dots~, (1, \dots,  1) 
		& 15  \\
\hline
8 & (8), (7, 1), (6, 2), (6, 1, 1), (5, 3), %(5, 2, 1), (5, 1, 1, 1),
	% (4, 4), 
	\dots~, (1, 1, 1, 1, 1, 1, 1, 1) 
		& 22  \\
\hline
\hline
~& ~& 60\\
\hline
\end{tabular}
\end{center} 
\end{table} 
%---------------------
% EXAMPLE
\begin{example}
Let us rewrite the network $2A \lra A+B \lra A+C$ as two separate reversible reactions:
	\begin{align} \label{eq:net_as_2_rxns}
	2A \lra A+B \quad \quad A+B \lra A+C~.
	\end{align}
Counting the number of times each species appears (where we take into consideration the stoichiometric coefficients), we see that species $A$ appears $5=2+1+1+1$ times, $B$ appears $2=1+1$ times, and $C$ appears 1 time.  Definition~\ref{def:TM} will say that the `total molecularities' of species $A$, $B$, and $C$ are, respectively, 5, 2, and 1.  In addition, the `total molecularity partition' of network~\eqref{eq:net_as_2_rxns} will be $(5,2,1)$, which is a partition of the integer $8=5+2+1$.  Similarly, the total molecularity partition of the network $2A \lra A+B \lra A$ is $(5,2)$, a partition of 7.  
\end{example}
The definition of total molecularity first appeared in \cite{Simplifying}.
% DEFINE TM
\begin{definition} \label{def:TM}
	\begin{enumerate}
	\item For a reversible network, the {\em total molecularity} of species $X_j$ refers to the sum over all pairs of reversible reactions of the sum of the stoichiometric coefficients of $X_j$ in the reactant and in the product:
	\begin{align*}
	\op{TM} (X_j) \quad := \quad \sum_{\sum_{i=1}^s a_i X_i \lra \sum_{i=1}^s b_i X_i ~\in~ \RR'} a_j+b_j~,
	\end{align*}
where $\RR'$ denotes all pairs of reversible reactions $\sum_{i=1}^s a_i X_i \lra \sum_{i=1}^s b_i X_i$.
	\item For a reversible network, the {\em total molecularity partition} is the partition defined by the multiset of total molecularities of all species: 
	\begin{align*}
	\{\op{TM}(X_1), ~\op{TM}(X_2),~\dots~,~ \op{TM}(X_{|\SSS|}) \}~.
	\end{align*}
	\end{enumerate} 
\end{definition}

Note that for reversible bimolecular two-reaction networks, the total molecularity partition is of an integer $m \in \{ 4,5,6,7,8\}$.  %The partitions of each such $m$ are listed (partially) in Table~\ref{table:ptn}. 
 
%------------------
% SUBSECTION: Algorithm
%------------------
\subsection{Algorithm for enumerating networks} \label{sec:algor}
We now present the algorithm we used for enumerating reversible bimolecular two-reaction networks.
	\begin{algorithm}[Algorithm for enumerating reversible bimolecular two-reaction networks] \label{algor:enum}
	~\\
	% STEP One
	{\bf Step One.}
	List partitions of $m=4,5,6,7,8$.\\
	% STEP Two
	{\bf Step Two.}
	For each partition $(m_1,m_2,\dots, m_n)$, list (with repeats) all reversible bimolecular two-reaction networks in which species $X_1$ has total molecularity $m_1$, species $X_2$ has total molecularity $m_2$, and so on.\\
	% STEP Three
	{\bf Step Three.}
	Remove networks that contain trivial reactions, networks that contain repeated reactions, and decoupled networks.\\
	% STEP Four
	{\bf Step Four.}
	Remove redundant networks: keep exactly one representative from each equivalence class of networks.  (Recall that two networks are equivalent if there exists a relabeling of the species that transforms the first network into the second network.)
	\end{algorithm}
As we see in Table~\ref{table:num_by_partition}, Algorithm~\ref{algor:enum} yields 386 reversible bimolecular two-reaction networks.  In Section~\ref{sec:class}, we determine which of the 386 CFSTRs admit multiple positive steady states.

Let us now elaborate on our implementations of Steps Two through Four of Algorithm~\ref{algor:enum}.  In order to list all reversible bimolecular two-reaction networks that have a given partition $(m_1,m_2,\dots, m_n)$ (Step Two), we made use of a psuedo-species $X_0$.  Namely, any network with partition $(m_1,m_2,\dots, m_n)$ arises from placing $(m-m_1-m_2-\cdots-m_n)$ copies of species $X_0$, $m_1$ copies of $X_1$, $m_2$ copies of $X_2$, and so on in the eight boxes in the following diagram:
	\begin{align*}
	\drawBox ~+~ \drawBox \quad \lra \quad \drawBox ~+~ \drawBox \quad \quad \quad \quad \quad \quad
	\drawBox ~+~ \drawBox \quad \lra \quad \drawBox ~+~ \drawBox 	%\framebox[1.1\height]{H} ~+~ 
	\end{align*}
For example, $\fbox{$A$}+\fbox{$A$} \lra \fbox{$A$}+\fbox{$B$} \quad \fbox{$A$}+\fbox{$B$} \lra \fbox{$X_0$}+\fbox{$A$}$ defines the network $2A \lra A+B \lra A$. Clearly, this procedure will yield all networks, but certain trivial networks (such as one with repeated reactions) will appear, and additionally each network will appear more than once.  Accordingly, trivial networks are removed in Step Three of Algorithm~\ref{algor:enum}, and Step Four keeps only one representative from each equivalence class of networks.  Step Four is the most computationally expensive part of our enumeration. {For each network remaining at the end of Step Three, we generated the equivalence class of networks obtained by performing a relabeling of the species. Two networks are equivalent if and only if they generate identical equivalence classes of networks. We removed extra copies of equivalent networks at the end of Step Four.}

%------------------
% SUBSECTION: Comparison with enumeration done by Deckard, Bergmann, and Sauro
%------------------
\subsection{The enumeration of small networks of Deckard, Bergmann, and Sauro} \label{sec:compare}
A related (and much larger: over 47 million) enumeration of small bimolecular networks was undertaken by Deckard, Bergmann, and Sauro~\cite{Deckard}.  Their work enumerated small networks by the number of {\em directed} reactions and by the number of species.  So, the network $2A \la A+B \la A$ falls in their list of networks containing two directed reactions and two species, and the network $2A \lra A+B \la A+C$ is a network containing three directed reactions and three species.  Also, their enumeration did not include seemingly unrealistic chemical reactions involving the zero complex (such as $0 \ra 2A$ or $0 \la A+B$) or reactions in which some species appears in both the reactant complex and product complex of a reaction (such as $A \ra A+B$ or $A \ra 2A$).  %[Add: how far they got - many networks! - and how their method of enumeration compares with ours: theirs via digraphs, etc.]
We remark that from this enumeration of networks by Deckard, Bergmann, and Sauro, the work of Pantea and Craciun sampled networks to compute the fraction that pass the Jacobian Criterion \cite[Figure 1]{Pantea_comput}.

%------------------
% SECTION: Classification
%------------------
\section{Classification of multistationary two-reaction CFSTRs} \label{sec:class}
The main result of this section is the following theorem, which completely answers Question~\ref{q:2RxnMSS}:
% MAIN CLASSIFICATION THEOREM
\begin{theorem} \label{thm:enum}
Of the {\bf 386} reversible, bimolecular, two-reaction \fd CFSTRs, exactly {\bf 35} admit multiple positive mass-action steady states.  Moreover, each of these 35 networks admits multiple {\em nondegenerate} positive steady states.  Furthermore, each such network contains a unique minimal multistationary subnetwork.  The poset (partially ordered set) of these 35 directed subnetworks, with respect to the embedded network relation, has {\bf 11} minimal elements, which are the bimolecular two-reaction CFSTR atoms of multistationarity.  %This poset is a poset (after adding $\hat 0$ and $\hat 1$).
\end{theorem}
An immediate corollary of Theorem~\ref{thm:enum} is the following:
% COROLLARY: 35 and 11
	\begin{corollary} \label{cor:class}
	\begin{enumerate}
	\item A bimolecular, two-reaction \fd CFSTR admits multiple nondegenerate positive mass-action steady states if and only if it contains as a sub-CFSTR one of the 35 minimal such subnetworks, which are displayed in Figure~\ref{fig:35}.
	\item A bimolecular, two-reaction \fd CFSTR admits multiple nondegenerate positive mass-action steady states if and only if it contains as an embedded network one of the 11 CFSTR atoms which are marked in bold/red in Figure~\ref{fig:35}.
	\item If a \fd CFSTR (not necessarily bimolecular and having any number of reactions) $G$ contains one of the 35 minimal CFSTRs mentioned above as a sub-CFSTR or contains one of the 11 atoms as an embedded network, then $G$ admits multiple nondegenerate positive mass-action steady states.
	\end{enumerate}
	\end{corollary}
\noindent
Note that part 3 of Corollary~\ref{cor:class} makes use of Corollaries~\ref{cor:liftCFSTR} and~\ref{cor:liftEmbCFSTR}.
	
% EXAMPLE
\begin{example}
Among the 35 reversible CFSTRs in Theorem~\ref{thm:enum} that admit multiple steady states, one is the network~\eqref{eq:running} which we first saw in Example~\ref{ex:running}: it arises from the network ${ 2A} \leftrightarrow { A+B} \leftrightarrow { A+C}$.  The 
unique minimal multistationary sub-CFSTR is the directed subnetwork~\eqref{eq:running_sub} obtained by removing two reactions: ${ 2A} \leftarrow { A+B} \leftarrow { A+C}$.  Finally, there is a multistationary embedded CFSTR~\eqref{eq:running_emb} obtained by removing species $C$, namely, the CFSTR arising from ${ 2A} \leftarrow { A+B} \leftarrow { A}$, that is one of the 11 atoms.
% minimal multistationary networks with respect to the `embedded relation'.   
 In other words, no further embedded CFSTR is multistationarity.  
 The directed subnetwork~\eqref{eq:running_sub} and the embedded network~\eqref{eq:running_emb} appear in the lower left of Figure~\ref{fig:35}.
\end{example}
Sections~\ref{subsec:rule_out} through~\ref{sec:irr} provide the proof for Theorem~\ref{thm:enum}.  

%----------
%:Table - counts number of networks by partition
\begin{table}
\begin{center}
\begin{tabular}{ |c|c|c|c|c|c|c|}
\hline
~ 		& Number of 		& Total 	& \#	&  \# 	& \# \\
~		& {\em reversible} bimolecular two-reaction networks	
					& \# of	& with		&  that fail  	& with \\	
$m$	 	& by partition of $m$	& networks	& TM $>2$	&  Jac. Crit. 	& MSS \\
\hline
4  & (0,2,{\bf 2,5,3}) & 12 & 2 & 0 & 0 \\
\hline
5  & (1,4,7,8,{\bf 10,9,2}) & 41 & 20 & 8 & 1\\
\hline
6  & (0,3,6,9,7,23,12,{\bf 9,23,12,3}) & 107 & 60 & 31 & 5 \\
\hline
7  & (0,1,3,4,5,13,7,9,13,26,8,{\bf 12,15,7,1}) & 124 & 89 & 55 & 15 \\
\hline
8  & (0,0,0,1,1,3,2,1,5,4,9,4,7,8,13,12,3,{\bf 5,11,9,3,1}) & 102 & 73 & 48 & 14 \\
\hline
\hline
~&~& 386  & 244 & 142 & {\bf 35} \\
\hline
\end{tabular}
\caption{Here we list the number of {\em reversible} bimolecular two-reaction CFSTRs by partition.  The order of partitions is the lexicographic order (as in Table~\ref{table:ptn}).  In bold are the 142 networks for which multistationarity is ruled out because each part of the corresponding partitions (the total molecularity of a species, denoted by ``TM'') is no more than two~\cite{Simplifying}.  Of the remaining 244 networks, an additional 102 networks pass the Jacobian Criterion (those with total molecularity at most two also pass the Jacobian Criterion).  For the remaining 142 networks, the CRN Toolbox~\cite{Toolbox} determined that precisely 35 admit multiple positive mass-action steady states (``MSS'').
\label{table:num_by_partition} }
\end{center} 
\end{table} 

%------------------
% SUBSECTION: Ruling out MSS by the JC
%------------------
\subsection{Ruling out multistationarity by the Jacobian Criterion} \label{subsec:rule_out}
Recall from~\cite{ME1,ME_entrapped,ME2,ME3} that the Jacobian Criterion is a method for ruling out multistationarity.  A CFSTR is said to pass the Jacobian Criterion if all terms in the determinant expansion of the Jacobian matrix of its mass-action differential equations~\eqref{eq:main} have the same sign.  Craciun and Feinberg proved that if a CFSTR passes the Jacobian Criterion, then it does not admit multiple positive steady states.  In earlier work, the current authors proved that if the total molecularities of all species are at most two, then the CFSTR passes the Jacobian Criterion~\cite{Simplifying}.  Accordingly, any two-reaction networks that arise from the 19 partitions (of 4, 5, 6, 7, or 8) in which all parts are at most two automatically pass the Jacobian Criterion; these 142 networks are marked in bold in Table~\ref{table:num_by_partition}.  Of the remaining $244=386-142$ networks, an additional 102 networks pass the Jacobian Criterion. % (including those with total molecularity at most two).

%We saw in Theorem~\ref{thm} that there are 386 reversible, two-reaction, bimolecular CFSTRs and of those, 142 have a maximum total molecularity of two or fewer, and thus immediately pass the Jacobian Criterion. 

%------------------------
% SUBSECTION: Apply CRNT to reversible networks that fail JC
%------------------------
\subsection{Applying the CRN Toolbox to classify reversible two-reaction networks}
For the remaining 142 reversible networks that do not pass the Jacobian Criterion, we applied the CRN Toolbox~\cite{Toolbox}.  This was performed in an automated fashion by using AutoIt code~\cite{AutoIt} provided by Dan Siegal-Gaskins.  We find that exactly 35 admit multiple positive mass-action steady states and the remaining 107 do not.  For each of the 35 multistationary CFSTRs, the Toolbox gave an instance of rate constants, two positive steady state values, and the corresponding eigenvalues.  In all cases but one, the nondegeneracy of these steady states was evident from the eigenvalues.  In the remaining case, in which one steady state was degenerate, we found `by hand' another instance of multistationarity in which two nondegenerate steady states exist. %Add more details: BB17 had one eigenvalue = 0

For the remaining 107 networks, the CRN Toolbox concluded that they do not admit multiple steady states.  A portion of a report produced by the Toolbox for such a network follows:
\begin{small}
\begin{verbatim}
Taken with mass action kinetics, the network CANNOT admit multiple positive 
steady states or a degenerate positive steady state NO MATTER WHAT (POSITIVE) 
VALUES THE RATE CONSTANTS MIGHT HAVE.
\end{verbatim}
\end{small}
The theoretical underpinning of the Toolbox consists of the Deficiency, Advanced Deficiency, and Higher Deficiency Theories developed by Ellison, Feinberg, Horn, Jackson, and Ji \cite{EllisonThesis,FeinDefZeroOne,HornJackson72,JiThesis}. 

%------------------------
% SUBSECTION: Apply CRNT to irreversible subnetworks
%------------------------
\subsection{Classifying irreversible two-reaction networks} \label{sec:irr}
\begin{figure}[h!]
\begin{center}
\includegraphics[angle=0,scale=1]{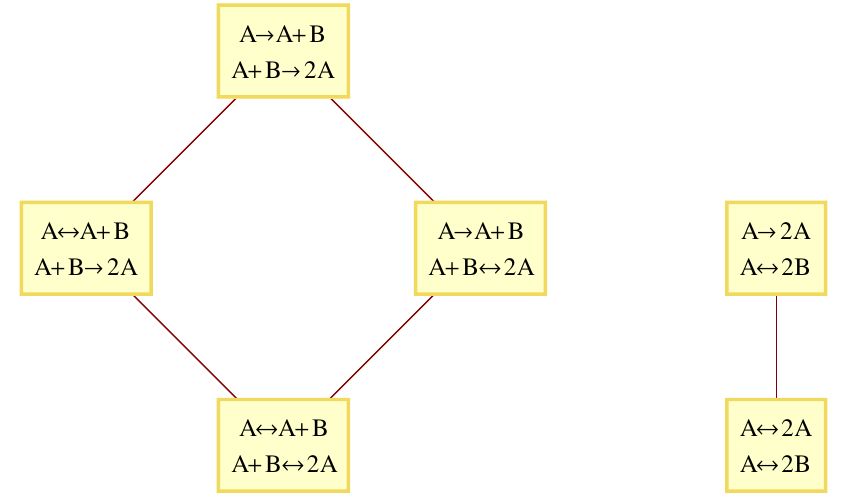} %formally called `Subnet1'
\end{center}
\caption{ Here we display two of the 35 multistationary bimolecular reversible two-reaction CFSTRs (all inflow and outflow reactions are implied), together with all their respective irreversible multistationary sub-CFSTRs (subnetworks).  The poset relation depicted is the subnetwork relation.  At left, the reversible CFSTR defined by reactions $A \lra A+B \lra 2A$ has three multistationary sub-CFSTRs, which are displayed above.  Similarly, the example on the right, defined by $2A \lra A \lra 2B$, has only one multistationary sub-CFSTR.  More generally, each of the 35 such reversible networks admits a unique minimal multistationary sub-CFSTR.  These minimal subnetworks fall into two classes: 29 of them have the form of the example displayed on the left (the minimal network has two directed reactions), and the remaining 6 have the form of the example on the right (the minimal network has one reversible reaction and one directed reaction).  These 35 minimal sub-CFSTRs appear in Figure~\ref{fig:35}. \label{fig:subnet}}
\end{figure}

Next, we consider the irreversible versions of the reversible two-reaction networks studied.  That is, we are interested in the networks obtained from the 386 reversible networks by making one or both of the non-flow reactions irreversible (each reversible reaction can be made irreversible in two ways).  So, each reversible network has 8 relevant subnetworks.  Recall that smaller sub-CFSTRs, those containing only one directed non-flow reaction or one pair of reversible non-flow reactions, were already analyzed in Theorem~\ref{thm:1rxn}, and the bimolecular hypothesis ensures that none are multistationary in the setting here.  By Theorem~\ref{thm:subnetwork}, only subnetworks of one of the 35 multistationary reversible networks can be multistationary.  Therefore, we must examine only $35*8=280$ such networks.  We again applied the Toolbox~\cite{Toolbox}.  We found that each of the 35 reversible CFSTRs has a unique minimal sub-CFSTR $N_i$ that admits multiple positive steady states.  Of these 35 subnetworks $N_i$, 29 of them have two directed non-flow reactions, while the remaining 6 have non-flow reactions that consist of 1 reversible reaction and 1 directed reaction.  Examples of both types appear in Figure~\ref{fig:subnet}.  {\em Thus, a bimolecular two-reaction (possibly irreversible) CFSTR admits multiple positive steady states if and only if one of these 35 minimal networks is a subnetwork} (part 1 of Corollary~\ref{cor:class}).

%-------------------------------------
% FIGURE OF 2 EXAMPLES of minimal subnetworks:

%-------------------------------------
% FIGURE OF 35 MINIMAL SUBNETWORKS:
% Poset with respect to `embedded' relation
\begin{figure}[p]
\begin{center}
\includegraphics[angle=0,scale=0.39]{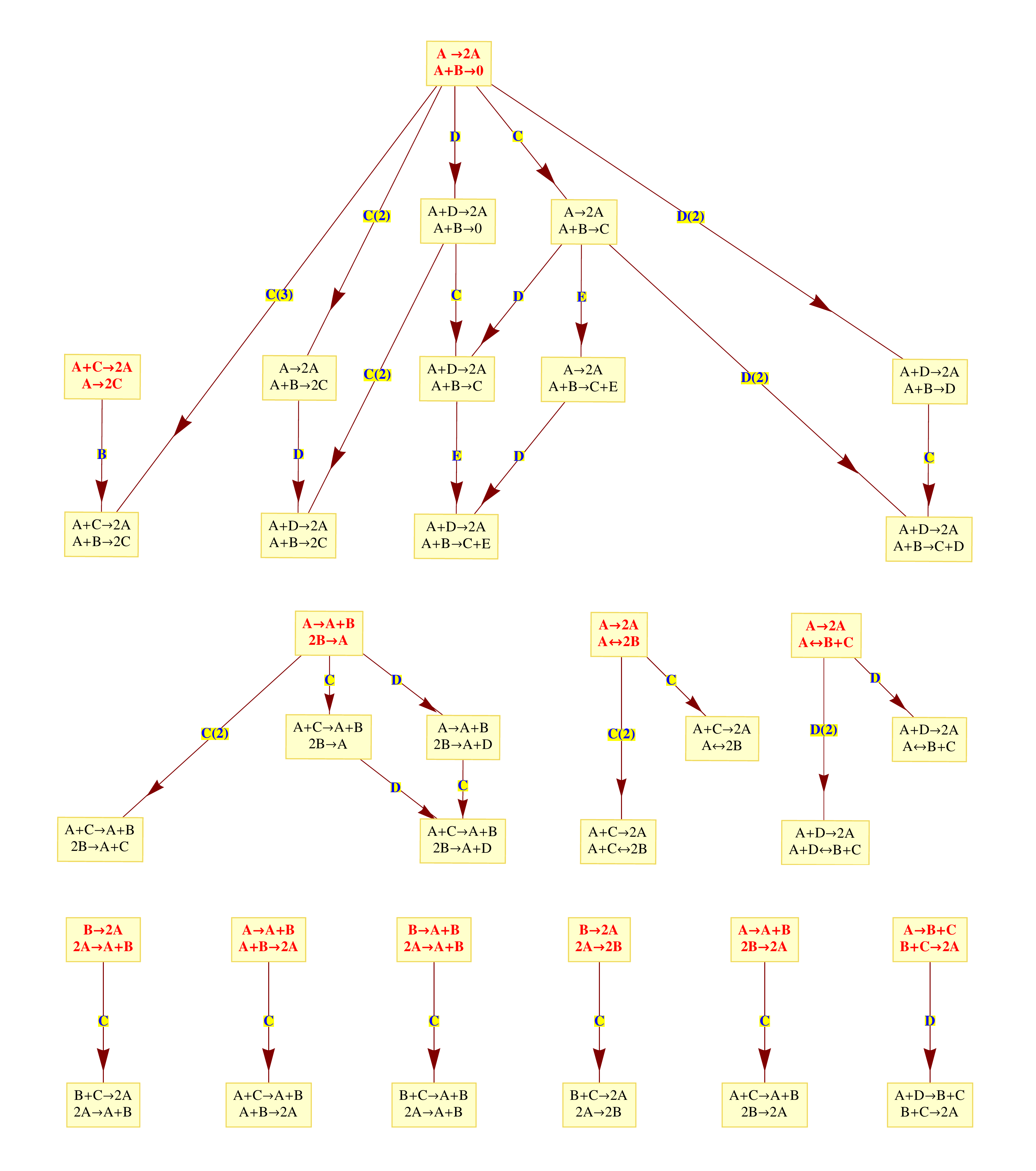} % formally called `LatticeAlt1'
\end{center}
\caption{Here we display the 35 multistationary bimolecular two-reaction CFSTRs that are minimal with respect to the subnetwork relation.  The poset relation depicted here is that of embedded networks: an arrow points from a network $N$ to a network $G$ if $N$ is an embedded network of $G$.  In addition, each such edge is labeled by the species that is removed to obtain $N$ from $G$; for example, $C(2)$ denotes that $G$ contains two molecules of species $C$, and these two are removed from $G$ to obtain $N$.  Two networks in the poset are displayed with the same height if they contain the same number of molecules.  The 11 CFSTR atoms of multistationarity are marked in bold/red; they are the networks that have only outgoing edges in the figure (at the tops of each component of the poset).  All three figures in this work were created in {\tt Mathematica}. \label{fig:35}}
\end{figure}

Finally, we examined the poset obtained from the $N_i$ with respect to the relation of `embedded networks' which is displayed in Figure~\ref{fig:35}.  This poset has 11 minimal elements, which are the bimolecular two-reaction CFSTR atoms of multistationarity.   {\em It follows that a bimolecular two-reaction (possibly irreversible) CFSTR admits multiple positive steady states if and only if it contains one of these 11 atoms as an embedded network} (part 2 of Corollary~\ref{cor:class}).

We end by noting that prohibitively many bimolecular {\em three}-reaction networks exist, so currently there is no classification of those CFSTR atoms. % of multistationarity

%[There are 3510 ($386 \times 9 + 12 \times 3$) bimolecular reaction networks containing at most two reversible or directed reactions.  *Could this include overcounting??]

%--------------------
\subsection*{Acknowledgements}
This project was initiated by Badal Joshi at a Mathematical Biosciences Institute (MBI) summer workshop under the guidance of Gheorghe Craciun. Badal Joshi was partially supported by a National Science Foundation grant (EF-1038593).  Anne Shiu was supported by a National Science Foundation postdoctoral fellowship (DMS-1004380).
The authors thank Dan Siegal-Gaskins for sharing AutoIt code which enabled the automated analysis of networks by the CRN Toolbox.

%--------------
% Bibliography
%--------------
\bibliographystyle{amsplain}
\bibliography{multistationarity}
%--------------
% END OF DOCUMENT
%--------------
\end{document}